\documentclass[draft]{amsart}
\usepackage{amssymb,amsmath,latexsym,epsfig,euscript, amsthm}

\def\boxit#1#2{\setbox1=\hbox{\kern#1{#2}\kern#1}%
\dimen1=\ht1 \advance\dimen1 by #1 \dimen2=\dp1 \advance\dimen2 by
#1
\setbox1=\hbox{\vrule height\dimen1 depth\dimen2\box1\vrule}%
\setbox1=\vbox{\hrule\box1\hrule}%
\advance\dimen1 by .4pt \ht1=\dimen1 \advance\dimen2 by .4pt
\dp1=\dimen2 \box1\relax}

\def\R{\mathbb{R}}

\def\Z{\mathbb{Z}}
\def\C{\mathbb{C}}

\def\N{\mathbb{N}}
\def\F{\mathbb{F}}
\def\CC{\mathcal{C}}

\def\AA{\mathcal{A}}
\def\EE{\mathcal{E}}

\def\div{\, |\,}
\def\ndiv{\nmid}

\def\Ima{{\rm Im}}
\def\Ker{{\rm Ker}}
\def\Rank{{\rm Rank}\,}

\def\set#1{{\left\{{\def\st{\;:\;}#1}\right\}}}

\theoremstyle{definition}
\newtheorem{defn}{Definition}[section]
\newtheorem{notn}[defn]{Notation}
\newtheorem{rem}[defn]{Remark}
\newtheorem{heur}[defn]{Heuristic}

\theoremstyle{plain}
\newtheorem{lem}[defn]{Lemma}

\newtheorem{prop}[defn]{Proposition}
\newtheorem{teo}[defn]{Theorem}
\newtheorem{alg}[defn]{Algorithm}
\newtheorem{cor}[defn]{Corollary}
\newtheorem{obs}[defn]{Observation}
\newtheorem{exmpl}[defn]{Example}

\newenvironment{undef}[1]%
           {\vspace{3.3mm}
           \noindent{\bf #1}\it}%
           {\vspace{3.3mm}}

\begin{document}

\title{Newton-Hensel interpolation lifting}

\author{Mart\'\i n Avenda\~no }
\address{Departamento de Matem\'atica, Facultad de Ciencias Exactas y
Naturales, Universidad de Buenos Aires, Pabell\'on I Ciudad
Universitaria 1428 Buenos Aires, Argentina}
\email{mavendar@yahoo.com.ar}
\thanks{The first and third author are partially supported by CONICET fellowships.}
\author{Teresa Krick}
\address{Departamento de Matem\'atica, Facultad de Ciencias Exactas y
Naturales, Universidad de Buenos Aires, Pabell\'on I Ciudad
Universitaria 1428 Buenos Aires, Argentina}
\email{krick@dm.uba.ar}
\thanks{The second author is supported by research grants UBACYT
X-112 and CONICET PIP 2461/01, Argentina.}
\author{Ariel Pacetti}
\address{Departamento de Matem\'atica, Facultad de Ciencias Exactas y
Naturales, Universidad de Buenos Aires, Pabell\'on I Ciudad
Universitaria 1428 Buenos Aires, Argentina}
\email{apacetti@dm.uba.ar}

\begin{abstract} The main result of this paper is a new version of
Newton-Hensel lifting
that relates to interpolation questions. It
 allows one to lift polynomials in $\Z[x]$
from information modulo a prime number $p\ne 2$ to a power $p^k$ for
any $k$, and its originality  is that it  is  a mixed version that
 not only lifts the coefficients of the polynomial but also its
 exponents. We show that this result corresponds exactly to a
Newton-Hensel lifting of a system of $2t$ generalized equations
in $2t$ unknowns in the ring of $p$-adic  integers $\Z_p$. Finally we apply our
results to sparse polynomial interpolation in $\Z[x]$.
\end{abstract}
\maketitle

\typeout{Introduccion}

\section*{Introduction}


\vspace{2mm}

Quoting \cite[p.10]{BCSS}, Newton's method has been considered
\textit{``the search algorithm'' sine qua non\/} of numerical
analysis and scientific computation. Since its origins probably by
Vi\`ete around 1580, its description by Newton in 1669,
simplification by Raphson in 1690 and actual formulation by
Simpson in 1740, Newton's method has been widely studied, applied
and generalized. For instance we mention here the crucial
development by S. Smale of the alpha theory, emphasizing
conditions on the input approximate zero (i.e., where can we start
the Newton's method) instead of hypotheses with regard to
estimates on the unknown zero, and of the gamma theory that
estimates the size of a ball of approximate zeros around the zero
\cite{Smale86} and \cite[Ch.~8 and 14]{BCSS}. Another important
issue is the search for algorithms for finding approximate zeros
in \cite{ShSm94, CuSm99} (the search of a polynomial time uniform
algorithm for such an approximate zero is one of the mathematical
problems for next century proposed in \cite{Smale98}), and the
generalizations of Newton's method for over-determined systems,
for instance in \cite{DeSh00}, and recently for systems with
multiplicities in \cite{Lecerf02}.

\smallskip
The non-archimedean counterpart of Newton's method, introduced by
Hensel  around 1900,  is the basis of  the construction  of the
$p$-adic
 numbers and their applications as in the local-global Hasse-Minkowski principle for quadratic
 forms. Since then,
 ``Newton-Hensel lifting" in its  non-archimedean versions is primordially
 present
 in exact symbolic computation: for example  in univariate rational polynomial
 factorization in \cite{Zassenhaus69} and the breakthrough LLL-polynomial time factorization
 algorithm
 \cite{LeLeLo82}, in multivariate polynomial
 factorization \cite{ChGr82,Chistov84,Grigoriev84,Kaltofen85}. Also for multivariate
 polynomial systems solving  in the Gr\"obner basis
 setting in \cite{Trinks85,Winkler88} and in the primitive element setting in
 \cite{ChGr83,Chistov84,Grigoriev84} and in
 \cite{GiHeMoMoPa98,GHHMMP97, HKPSW00, JKSS04}.

\smallskip
The main result of this paper is a new version of Newton-Hensel lifting
that relates to interpolation questions. It
 allows to lift polynomials in $\Z[x]$
from information modulo a prime number $p\ne 2$ to a power $p^k$ for
any $k$, and its originality  is that it  is  a mixed version that
 not only lifts the coefficients of the polynomial but also its {\em
 exponents}.

\begin{undef}{Theorem 1.} Let $p$ be an odd  prime number, and  $t\in \N$, $1\le t\le
{p-1}$.

Set  $y_1,\dots, y_{2t}\in \Z$. Let
$f:=\sum_{j=1}^ta_jx^{\alpha_j}\in\Z[x]$ and $x_1,\dots,x_{2t}\in
\Z$ satisfy that

\begin{itemize}
\item $f(x_i)\equiv y_i \pmod p$, $1\le i\le 2t$,
\item $\det\,\left(
\begin{array}{cccccc}x_1^{\alpha_1} & \dots &
x_1^{\alpha_t}&a_1\,e_p(x_1)x_1^{\alpha_1}& \dots &
a_t\,e_p(x_1)x_1^{\alpha_t}\\
 \vdots  & & \vdots  & \vdots & & \vdots \\
x_{2t}^{\alpha_1} &  \dots &
x_{2t}^{\alpha_t}&a_1\,e_p(x_{2t})x_{2t}^{\alpha_1}&
 \dots & a_t\,e_p(x_{2t})x_{2t}^{\alpha_t}\end{array}\right)
\not\equiv 0 \pmod p$, \\[2mm]
where $e_p:\Z \to \Z/p\Z$ is the map defined by $x^{p-1}\equiv 1+p
\,e_p(x) \pmod{p^2}$ for $p\ndiv x$ and $e_p(x)=0$ for $p\div x$
(see Definition \ref{morphism} and Notation \ref{e-notation}).
\end{itemize}

Then for every $k\in \N_0$  there exists $f_k:=\sum_{j=1}^ta_{k,
 j}x^{\alpha_{k, j}}\in \Z[x]$, that satisfies simultaneously:
\begin{itemize}
\item $ f_k(x_i) \equiv  y_i \pmod{p^{2^k}}$ for $ 1\le i\le 2t$,

\item $a_{k, j}\equiv a_j \pmod p$ and $\alpha_{k ,j}\equiv
\alpha_j\pmod{(p-1)}$, $1\le j\le t$.
\end{itemize}

\smallskip
Furthermore, if $f_k:=\sum_{j=1}^ta_{k, j}x^{\alpha_{k, j}},
g_k:=\sum_{j=1}^tb_{k ,j}x^{\beta_{k ,j}}\in\Z[x]$ are two such
polynomials, then $$ a_{k, j}\equiv b_{k ,j} \pmod{p^{2^k}} \
\mbox{ and } \ \alpha_{k, j}\equiv \beta_{k,
j}\pmod{\varphi(p^{2^k})}, \ 1\le j\le t,$$ where $\varphi$
denotes the Euler map.
\end{undef}

 In
fact we show in this paper that this result corresponds exactly to  a
Newton-Hensel lifting of a system of $2t$ generalized equations
in $2t$ unknowns, where the unknowns are the $t$ coefficients of $f$ in
the $p$-adic integers $\Z_p$ and the $t$ exponents of $f$ in  some
suitable set $\EE_p$ (see Definition \ref{Ep}), and where the condition that the defined
 determinant does not vanish modulo $p$ is the corresponding classical
condition on the Jacobian determinant of the system: the map $e_p$
plays the role of a logarithmic function that  enables us to lower
the  exponents to the floor level.

\smallskip
 For this purpose we introduce the ring
$\EE_p$ of ``allowed exponents" (whose additive group is
isomorphic to the $p$-adic unit group $\Z_p^\times$), where ``allowed" means that
  for $x\in \Z_p^\times$ and $\alpha\in \EE_p$, $x^\alpha \in
  \Z_p^\times$, and
  we study systems of generalized polynomial
expressions in $\Z_p$, where the variables belong to $\Z_p^\times$
and  the exponents belong to  $\EE_p$.

\smallskip
  Here, among all the equivalent descriptions of $\Z_p$ we adopt  the following one
  that we consider more
  suitable for our formulations:
$$
\mathbb{Z}_{p}=\left\{ (a_{k})_{k\in \N}\in\mathbb{Z}^{\mathbb{N}}\;
:\; a_{k+1}\equiv a_{k}\pmod{p^{k}}\ \forall\, k\in \N\,\right\}
/\sim ,$$ where $
 (a_{k})_{k\in \N}\sim(b_{k})_{k\in \N}\ \Leftrightarrow \ a_{k}\equiv
b_{k}\pmod{p^{k}},\;\forall\, k\in \N$. Similarly, we have
$$
\EE_{p}=\left\{ (\alpha_{k})_{k\in \N}\in\mathbb{Z}^{\mathbb{N}}\;
:\; \alpha_{k+1}\equiv \alpha_{k}\pmod{\varphi(p^{k})}\ \forall\,
k\in \N\,\right\} /{\approx} ,$$ where $
 (\alpha_{k})_{k\in \N}\approx (\beta_{k})_{k\in \N}\ \Leftrightarrow \ \alpha_{k}\equiv
\beta_{k}\pmod{\varphi(p^{k})},\;\forall\, k\in \N$. (In the sequel
$\overline{(a_k)}_k$ or $\overline{(\alpha_k)}_k$ denote the class
in the corresponding ring.)

\smallskip
 We consider   systems of equations where the unknowns are the variables in $\Z_p^\times$
 or the exponents in $\EE_p$, switching from one formulation to the other by a
 logarithmic
 argument,
 and obtain in Propositions \ref{henselsyspolyequation} and \ref{henselsysexpequation}
 below the generalizations of the following Newton-Hensel univariate lifting statements:

 \begin{undef} {Proposition 2.}
 Set  $t\in \N, \,y,\, a_j\in \Z_p$
and $\alpha_j:=\overline{(\alpha_{jk})}_k\in \EE_p$, $1\le j\le t$.
Let $x_1\in  \Z$, $p\ndiv x_1$,  be such that
\begin{itemize}
\item $f(x_1):=\sum_{j=1}^ta_j\,x_1^{\alpha_j}\equiv y \pmod p$,
\item $\sum_{j=1}^t a_{j}\,\alpha_{j2}\,x_1^{\alpha_{j}-1}\not\equiv 0 \pmod p$,
\end{itemize}
then there exists a unique $x=\overline{(x_1,x_2,\dots )}\in \Z_p^\times$ such that
$f(x)=y$ in $\Z_p$.
\end{undef}

\begin{undef} {Proposition 3.}
 Set  $t\in \N,\, y, \,a_j\in \Z_p$
and $x_j:=\overline{(x_{jk})}_k\in \Z_p^\times$, $1\le j\le t$.
Let $\alpha_1\in \Z$ be such that
\begin{itemize}
\item $g(\alpha_1):=\sum_{j=1}^ta_j\,x_j^{\alpha_1}\equiv y \pmod p$,
\item $\sum_{j=1}^t a_{j}\,e_p(x_{j2})\,x_{j}^{\alpha_{1}}\not\equiv 0 \pmod p$,

where $e_p:\Z \to \Z/p\Z$ is the map defined by $x^{p-1}\equiv 1+p
\,e_p(x)\pmod{p^2}$ for $p\ndiv x$ and $e_p(x)=0$ for $p\div x$,
\end{itemize}
then there exists a unique $\alpha=\overline{(\alpha_1,\alpha_2,\dots)}\in \EE_p$ such that
 $g(\alpha)=y$ in $\Z_p$.
 \end{undef}

We observe that the  map $e_p$ used in the hypotheses of Theorem 1
and Proposition 3
 is a group
homomorphism $(\Z/p^2\Z)^\times \to \Z/p\Z$ that appears widely
in the literature, even
 if not explicitly stated in the
way we do  in Proposition \ref{e-morphism} below. It arises
naturally when looking for generators of the cyclic multiplicative
group $\Z/p^k\Z$ for $p$ an odd prime number
\cite[Th.10.6]{Apostol76} or  considering the filtration of the
group of $p$-adic units $\Z_p^\times$ \cite[Ch.II \S
3.1]{Serre70}. As explained in \cite[p.101]{Carmichael14} the
properties $e(2)= 0$ and $e(3)=0$ are also intimately related to
special cases of Fermat's last theorem. It seems that  at least
until 1999, it was still unknown if there were infinitely many
prime numbers $p$ such that $e_p(2)\ne 0$ without assuming the ABC
conjecture \cite[p.8]{EsMu99}. Finally, let us mention that
numerical experiments we made suggest that the behavior of the map
$e_p(x)/p-1$ for a fixed integer $x$ and a variable prime number
$p$ seems to follow a uniform distribution in the $[0,1]$
interval.

\bigskip
Our initial motivation and a possible useful application for the
Newton-Hensel lifting result presented in Theorem 1 above was
 the search of an efficient interpolation algorithm for integer
 univariate
 polynomials, where the number of interpolation points depends on the number of
 non-zero terms and not on the degree.
A  polynomial is called {\em $t$-sparse} if it has at most $t$
non-zero terms. The problem of interpolating  a $t$-sparse
polynomial from its values in a list of specific interpolation
points  where the number of these points does not depend on the
degree but on $t$ is called ``sparse interpolation". It  received a
lot of attention around 1990 and again recently, for instance in
\cite{BeTi,Zip90,KaLa88, KLW90, KLL00,Lee01,KaLe03}, \cite{BoTi},
\cite{DrGr91,GKS91}, \cite{GrKaSi90,ClDrGrKa91} and \cite{GKS94}.

\smallskip

It is a well-known fact that, as a consequence of Descartes rule of
signs, a $t$-sparse polynomial $f\in \R[x]$ (one variable) has at
most $t-1$ distinct real positive roots. Therefore, any univariate
$t$-sparse polynomial in $\C[x]$ is uniquely determined by its value
in $2t$ different positive values in $\R$ (since for two such
polynomials, the difference $f-g$ of their real parts (or their
imaginary parts) is a $2t$-sparse polynomial which has at most
$2t-1$ different real positive roots). In \cite{BeTi},  M. Ben-Or
and P. Tiwari produced a beautiful deterministic algorithm that
recovers such a $t$-sparse polynomial $f\in \C[x]$ from its value in
the $2t$ interpolation points
$$x_1:= 1,\,
x_2:=a,\,x_3:=a^2,\dots,x_{2t}:=a^{2t-1},$$ where $a$ is not a root
of unity of small order. They also  raised the problem of producing
an algorithm that interpolates a $t$-sparse polynomial in $\C[x]$
from $2t$ arbitrary different real positive values, to emulate in
some sense Lagrange or Newton interpolation algorithms that do not
require specific interpolation input values, instead of imposing the
starting points as they do.

\smallskip
Although we are not able to answer this question in generality, we
produce in Proposition \ref{startingset} families
$\{x_1,\dots,x_{2t}\}$ of starting points where the non-vanishing
determinant hypothesis holds under a good reduction property of the
input $t$-sparse polynomial modulo the prime number $p$ (Proposition
\ref{pkequivalencia}). Therefore, applying Theorem 1 we obtain a
very fast algorithm for sparse interpolation of $t$-sparse
polynomials in $\Z[x]$ that
 reduce well modulo $p$ (Algorithm \ref{interpolation}). The algorithm does not require
 to know in advance the degree of the polynomial although it needs to know  its
  exact number  $t$ of non-zero terms.
  In order to
make this algorithm work for any   polynomial $f\in\Z[x]$,
 we still idealistically need a criterion to choose, in terms of the evaluation values $y_i=f(x_i)$,
  a (small) prime number  $p$
 such that $f$ reduces well modulo $p$. This would be
 the analog of the choice of the prime in the univariate rational
 polynomial
 factorization algorithms (the condition in this case is given by the
 non-vanishing of a discriminant modulo $p$), and in the archimedean
 setting, a (still unknown) criterion  for the choice of  an approximate zero.
 More realistically, we would
  at least need a satisfactory probabilistic
 argument for the choice of  such a prime in a
 given range (that we are still unable to produce).

 \bigskip
 The paper is organized as follows.

Section \ref{Newton-Hensel lifting} is mainly devoted to the proof
of Theorem 1. For this purpose we introduce a generalization of
the $e_p$ group homomorphism mentioned above (Definition
\ref{e-morphism}). Then we prove the theorem (Theorem
\ref{hensel}) and present an equivalence of the uniqueness
condition (Proposition \ref{pkequivalencia}). Finally we analyze
the existence of good starting sets $\{x_1,\dots, x_{2t}\}$ as
inputs of Theorem 1 (Definition \ref{defstartingset} and
Proposition \ref{startingset}).

In Section \ref{$p$-adic equations}, we focus on the Hensel lemma
character of our Theorem $1$ in the ring of $p$-adic integers
$\Z_p$. We introduce the set $\EE_p$ of allowed exponents
(Definition \ref{Ep}), the generalized polynomial equations and
their dual exponential equations (Observation \ref{dual}) and we
present the proofs of Propositions 1 and 2 above (Propositions
\ref{henselpolyequation} and \ref{henselexpequation}) and their
generalizations to systems of generalized polynomial and exponential
equations (Propositions \ref{henselsyspolyequation} and
\ref{henselsysexpequation}).

Finally Section \ref{Fewnomial interpolation} deals with the sparse
interpolation problem mentioned above, first focusing in univariate
$t$-sparse polynomials  with coefficients in  finite rings and then
in univariate   integer $t$-sparse polynomials.

 \bigskip
 {\bf Acknowledgments.} Teresa Krick thanks Wen-shin Lee for bringing
 her into the sparse interpolation subject and  Michael Singer
 for wonderful discussions. We are also grateful to  Felipe Cucker,
 Arieh Iserles,
   Erich Kaltofen and Mike Shub
  for their advice, and to the referees for their careful reading
  and clarifying comments.

\section{Newton-Hensel lifting}
\label{Newton-Hensel lifting}

This section is mainly devoted to the proof of   Theorem 1,
the Newton-Hensel interpolation lifting theorem stated in the introduction.

 \smallskip During
the paper $p$ denotes an odd prime number. Given an integer $\rho$
prime to $p$,  we denote by $o_{p^k}(\rho)$ the order of $\rho$ in
the multiplicative cyclic group $(\Z/p^k\Z)^\times$. We recall
that $\rho \in \Z$ is a primitive root modulo ${p^k}$ if its class
in $\Z/p^k\Z$ generates this cyclic group, that is if
$o_{p^k}(\rho)=\varphi(p^k)=(p-1)p^{k-1}$.

\smallskip
The crucial tool in this paper is a family of group homomorphisms
that relates the multiplicative group $(\Z/p^k\Z)^\times$ with the
additive group $\Z/p^{k-\ell}\Z$ for $\ell\le k\le 2\ell$. This morphism enables us to
linearize polynomial expressions.

\begin{defn}
Let $p$ be an odd prime number, $k,\ell$ positive integers with
$\ell \le k \le 2\ell$. Define the morphism $e_{p,k,\ell} : (\Z/p^k
\Z)^\times \rightarrow \Z/p^{k-\ell}\Z$ by
\[
e_{p,k,\ell}(x) := \frac{x^{\varphi(p^\ell)}-1}{p^\ell} \pmod {p^{k-\ell}}
\]
\label{morphism}
\end{defn}

\begin{prop} The map $e_{p,k,\ell}$ defined above is a group epimorphism.
\label{e-morphism}
\end{prop}

\begin{proof}
The map is clearly well-defined  since $x^{\varphi(p^\ell)}\equiv 1
+ p^{\ell}e_{p,k,\ell}(x) \pmod{p^{k}}$. Hence we must prove  it is
a group homomorphism. It is enough to prove:

\begin{itemize}
\item $e_{p,k,\ell}(1)=0$ (clear from the definition).
\item $e_{p,k,\ell}(x\, y ) = e_{p,k,\ell}(x)+e_{p,k,\ell}(y)$. Since
$x^{\varphi(p^{\ell})}\equiv 1+ p^{\ell}\,e_{p,k,\ell}(x) \pmod{p^{k}}$ and
$y^{\varphi(p^{\ell})} \equiv 1 + p^{\ell}\,e_{p,k,\ell}(y) \pmod{p^{k}}$, it follows
$(x\,y)^{\varphi(p^{\ell})} \equiv (1+p^\ell\,e_{p,k,\ell}(x))(1+p^{\ell}\,e_{p,k,\ell}(y))\equiv
1+p^{\ell}\,(e_{p,k,\ell}(x)+e_{p,k,\ell}(y)) \pmod{p^{k}}$ (here we use the condition $k
\le 2\ell$).
\end{itemize}

To see it is surjective we compute the order of its
kernel. $\Ker(e_{p,k,\ell})=\{x\in (\Z/p^{k}\Z)^\times :
x^{\varphi(p^{\ell})}\equiv 1 \pmod{p^{k}}\}= \{x\in
(\Z/p^{k}\Z)^\times : o_{p^{k}}(x)\div
\varphi(p^{\ell})=(p-1)p^{\ell-1}\}$. Then $|\Ker(e_{p,k,\ell})|=(p-1)p^{\ell-1}$,
$|\Ima(e_{p,k,\ell})|=p^{k-\ell}$ and $e_{p,k,\ell}$ is surjective.
\end{proof}

\begin{rem} If $\rho$ is a primitive root modulo $p^{k}$ then
$e_{p,k,\ell}(\rho)\ne 0$.
\end{rem}
\begin{notn} Abusing notation we will denote $e_{p,k,\ell}$ the map from
$\Z\rightarrow \Z/p^k\Z \rightarrow \Z/p^{k-\ell}\Z$ defined as in
Definition \ref{morphism} on elements prime to $p$ and by zero on
other elements. In the case $k=2$ we will denote $e_{p,2}(x) :=
e_{p,2,1}(x)$ (or just $e_2(x)$ if the prime is clear from the
context). \label{e-notation}
\end{notn}

\begin{defn} Let $p$ be an odd prime number, and $t\in \N$,  $1\le t\le p-1$.
Let $f=\sum_{j=1}^t a_{j}x^{\alpha_{j}}\in\Z[x]$.  For
$x_1,\dots,x_{2t}\in\Z$, we define:
\begin{equation}\label{deltap}
\Delta_p\, f (x_1,\dots,x_{2t}):= \det\,\left(
\begin{array}{cccccc}x_1^{\alpha_1} & \dots &
x_1^{\alpha_t}&a_1\,e_2(x_1)x_1^{\alpha_1}& \dots &
a_t\,e_2(x_1)x_1^{\alpha_t}\\
 \vdots  & & \vdots  & \vdots & & \vdots \\
x_{2t}^{\alpha_1} &  \dots &
x_{2t}^{\alpha_t}&a_1\,e_2(x_{2t})x_{2t}^{\alpha_1}&
 \dots & a_t\,e_2(x_{2t})x_{2t}^{\alpha_t}
\end{array}\right) \pmod p
\end{equation}

where $e_2(x_i):\Z \rightarrow \Z/p^2\Z \rightarrow \Z/p\Z$ is the
map of Notation \ref{e-notation}.
\end{defn}

Quantity (\ref{deltap}) plays the role of the usual jacobian in our
version of Newton-Hensel interpolation lifting. In the next section,
using the analysis in the $p$-adic context, the relation will become
clear. For that reason we refer to it as the {\em interpolating
pseudo-jacobian modulo $p$ of $f$ on $x_1,\dots,x_{2t}$}.

\begin{teo} \label{hensel}{\em (Newton-Hensel interpolation lifting)}
Let $p$ be an odd  prime number, and  $t\in \N$, $1\le t\le {p-1}$.
Set  $y_1,\dots, y_{2t}\in \Z$. Let
$f:=\sum_{j=1}^ta_jx^{\alpha_j}\in\Z[x]$ and $x_1,\dots,x_{2t}\in
\Z$ satisfy

\begin{itemize}
\item $f(x_i)\equiv y_i $ $\pmod p$, $1\le i\le 2t$, \item
$\Delta_p\, f(x_1,\dots,x_{2t}) \not\equiv 0$ $\pmod p$, where
$\Delta_p\,f$ is given by Formula (\ref{deltap}).
\end{itemize}

Then for every $k\in \N_0$  there exists
$f_k:=\sum_{j=1}^ta_{k,j}x^{\alpha_{k,j}}\in \Z[x]$, that satisfies
simultaneously:
\begin{itemize}
\item $ f_k(x_i) \equiv  y_i$ $\pmod {p^{2^k}}$ for $
1\le i\le 2t$,

\item $a_{k,j}\equiv a_j $ $\pmod p$ and $\alpha_{k,j}\equiv \alpha_j$
  $\pmod {(p-1)}$, $1\le j\le t$.
\end{itemize}

\smallskip
Furthermore, if $f_k:=\sum_{j=1}^ta_{k,j}x^{\alpha_{k,j}},
g_k:=\sum_{j=1}^tb_{k,j}x^{\beta_{k,j}}\in\Z[x]$ are two such
polynomials, then \begin{equation}\label{uniqueness} a_{k,j}\equiv
b_{k,j} \pmod{p^{2^k}} \ \mbox{ and } \ \alpha_{k,j}\equiv
\beta_{k,j}\pmod{\varphi(p^{2^k})}, \ 1\le j\le t. \end{equation}
\end{teo}

\smallskip
\begin{proof}
We define $f_0:=f$, which is clearly the only possible definition
for $f_0$. Assume now $f_k= \sum_{j=1}^{t} a_{k,j}x^{\alpha_{k,j}} $
is uniquely defined under Condition (\ref{uniqueness}), with $
f_k(x_i)\equiv y_i \pmod{p^{2^{k}}}$ for $1\le i\le 2t$, and
$a_{k,j}\equiv a_j \pmod p$, $\alpha_{k,j}\equiv \alpha_j \pmod
{(p-1)}$.

We look for $ f_{k+1}$ such that $f_{k+1}(x_i)\equiv y_i$
mod${p^{2^{k+1}}}$ for $1\le i\le 2t$. In particular
$f_{k+1}(x_i)\equiv y_i \pmod{p^{2^{k}}}$ for $1\le i\le 2t$, and
therefore $f_{k+1}$ is of the form

$$
 f_{k+1}:=  \sum_{j=1}^{t}(a_{k,j}+ p^{2^{k}}d_j) x^{(\alpha_{k,j}+
\varphi(p^{2^k})\delta_j)}
$$
 and
$d_j,\delta_j $ have to be determined such that $ f_{k+1}(x_i)\equiv
y_i\pmod{{p^{2^{k+1}}}}$ for $1\le i\le 2t$.

For notational simplicity, for the rest of the proof we set
$b_j:=a_{k,j}$ and $\beta_j:=\alpha_{k,j}$.

Since $x_i\not \equiv 0 \pmod p$, $x_i^{\varphi(p^{2^k})}\equiv 1
\pmod {p^{2^{k}}}$. Let $e_{2^{k+1}}(x_i), \ell_{2^{k+1}}(x_i)\in \{0,\dots, p^{2^{k}}-1\}$
be
such that
\begin{equation}\label{eielei}
x_i^{\varphi(p^{2^k})}\equiv 1 + p^{2^{k}}e_{2^{k+1}}(x_i)\pmod{p^{2^{k+1}}} \
\mbox{ and } \ y_i\equiv f_k(x_i)+p^{2^{k}}\ell_{2^{k+1}}(x_i)
\pmod{p^{2^{k+1}}}
\end{equation}
i.e. $e_{2^{k+1}}(x_i) = e_{p,2^{k+1},2^k}(x_i)$. Then

\begin{eqnarray*}
y_i&\equiv &   f_{k+1}(x_i)\pmod{p^{2^{k+1}}}\\
&\equiv & \sum_{j=1}^{t}(b_{j}+ p^{2^{k}}d_j ) x_i^{(\beta_j+
{\varphi(p^{2^k})}\delta_j)}
\pmod{p^{2^{k+1}}} \\
 &\equiv &\sum_{j=1}^{t}(b_j+ p^{2^{k}}d_j)
x_i^{\beta_j}(1+p^{2^k}e_{2^{k+1}}(x_i))^{\delta_j} \pmod{p^{2^{k+1}}}\\
&\equiv & \sum_{j=1}^{t}(b_{j}+ p^{2^{k}}d_j)
x_i^{\beta_j}(1+p^{2^k}e_{2^{k+1}}(x_i)\delta_j) \pmod{p^{2^{k+1}}}\\
&\equiv & \sum_{j=1}^{t}(b_{j}x_i^{\beta_j} +
p^{2^k}(x_i^{\beta_j}d_j +
b_je_{2^{k+1}}(x_i)x_i^{\beta_j}\delta_j
))\pmod{p^{2^{k+1}}}\\
&\equiv & f_k(x_i) + p^{2^{k}}\sum_{j=1}^t(x_i^{\beta_j}d_j +
b_je_{2^{k+1}}(x_i)x_i^{\beta_j}\delta_j)
\pmod{p^{2^{k+1}}}\\
&\equiv& y_i - p^{2^{k}} \ell_{2^{k+1}}(x_i) +
p^{2^{k}}\sum_{j=1}^t(x_i^{\beta_j}d_j +
b_je_{2^{k+1}}(x_i)x_i^{\beta_j}\delta_j )\pmod{p^{2^{k+1}}}.
\end{eqnarray*}
 Then
$$-p^{2^{k}} \ell_{2^{k+1}}(x_i) + p^{2^{k}}\sum_{j=1}^t(x_i^{\beta_j}d_j +
b_je_{2^{k+1}}(x_i)x_i^{\beta_j}\delta_j )\equiv 0
\pmod{p^{2^{k+1}}}$$ Dividing by $p^{2^k}$ we get
$$-\ell_{2^{k+1}}(x_i) + \sum_{j=1}^t(x_i^{\beta_j}d_j +
b_je_{2^{k+1}}(x_i)x_i^{\beta_j}\delta_j )\equiv 0 \pmod{p^{2^{k}}}\quad
\mbox{for} \ 1\le i\le 2t.$$

\smallskip
Thus, one has to solve modulo $p^{2^{k}}$ the linear system of
equations

\begin{equation}\label{sistema}
\begin{array}{cccc}  \left( \begin{array}{cccccc}x_1^{\beta_1} &
\dots & x_1^{\beta_t}&b_1e_{2^{k+1}}(x_1)x_1^{\beta_1}& \dots &b_te_{2^{k+1}}(x_1)x_1^{\beta_t}\\
 \vdots &  & \vdots & \vdots & & \vdots \\
x_{2t}^{\beta_1} &  \dots &
x_{2t}^{\beta_t}&b_1e_{2^{k+1}}(x_{2t})x_{2t}^{\beta_1}&
 \dots & b_te_{2^{k+1}}(x_{2t})x_{2t}^{\beta_t}
\end{array}\right)& \left(\begin{array}{c}d_1\\ \vdots\\ d_t
\\ \delta_1  \\ \vdots \\ \delta_t \end{array}\right) &=& \left( \begin{array}{c}\ell_{2^{k+1}}(x_1)\\
\vdots \\
\ell_{2^{k+1}}(x_{2t})\end{array}\right)
\end{array}
\end{equation}

Let $M_k$ denote the $2t$ square matrix on the left hand side of
System (\ref{sistema}), that is the matrix of the
$(k+1)$-iteration step of our construction. Our aim is to show
that this matrix is invertible modulo $p^{2^k}$, i.e.  that
$\det(M_k)\not\equiv 0 \pmod{p}$. As the next lemma shows all
matrices reduce to $M_0$ modulo $p$, then we can  restrict
ourselves to study $\det(M_0)$.

\begin{lem} With the same notation as above $M_k\equiv M_0\pmod p$
  term by term.
\end{lem}

\begin{proof}
By construction it is clear that $b_j \equiv a_j $ mod${p}$ and
$x_i^{\beta_j}\equiv x_i^{\alpha_j}$ $\pmod {p}$. \\[2mm]
Hence we are left to prove that for any $x \in \Z$ prime to $p$, $e_{p,2^{k+1},2^k}(x) \equiv
e_{p,2^{k+2},2^{k+1}}(x) \pmod{p^{2^k}}$. Let $e :=
e_{p,2^{k+1},2^k}(x)$ and $e' :=e_{p,2^{k+2},2^{k+1}}(x)$. Then
$x^{\varphi(p^{2^k})}\equiv 1+ p^{2^k}e
\pmod{p^{2^{k+1}}} $ and $x^{\varphi(p^{2^{k+1}})}\equiv 1
+p^{2^{k+1}}e' \pmod{p^{2^{k+2}}}$. Since
$\varphi(p^{2^{k+1}})=\varphi(p^{2^k})p^{2^k}$, we have:
\begin{eqnarray*}
x^{\varphi(p^{2^{k+1}})}  &= & (x^{\varphi(p^{2^k})})^{p^{2^k}} \\
&=&(1+ p^{2^k}(e+p^{2^k}r))^{p^{2^k}}\\
&=&1+ p^{2^k}p^{2^{k}}(e+p^{2^k}r)+ {p^{2^k}\choose
2}{(p^{2^k})^2}(e+p^{2^k}r)^2
+{p^{2^k}\choose 3}{(p^{2^k})^3}(e+p^{2^k}r)^3 +\cdots \\
&=&1+ p^{2^{k+1}}}(e+p^{2^k}r)+
{p^{2^k+2^{k+1}}\frac{p^{2^k}-1}{2}(e+p^{2^k}r)^2
+{p^{3\,2^k}}{p^{2^k}\choose 3}(e+p^{2^k}r)^3 +\cdots \\
&=&1+p^{2^{k+1}}\Big(e + p^{2^k}(r+\frac{p^{2^k}-1}{
2}(e+p^{2^k}r)^2+
{p^{2^k}\choose 3}(e+p^{2^k}r)^3+\cdots ) \Big)\\
&\equiv&1+p^{2^{k+1}}\Big(e + p^{2^k}r'\Big) \pmod{p^{2^{k+2}}},
\end{eqnarray*}
where $r':= r+ \dots$. Then $e' \equiv e + p^{2^k}r' \equiv e  $
mod${p^{2^k}}$.
\end{proof}

We continue with the proof of the theorem.  We can restrict to
compute $\det(M_0)$. By definition,
$$M_0=  \left( \begin{array}{cccccc}x_1^{\alpha_1} &
\dots & x_1^{\alpha_t}&a_1e_2(x_1)x_1^{\alpha_1}&
\dots &
a_te_2(x_1)x_1^{\alpha_t}\\
  \vdots & & \vdots  & \vdots & & \vdots \\
x_{2t}^{\alpha_1} &  \dots &
x_{2t}^{\alpha_t}&a_1e_2(x_{2t})x_{2t}^{\alpha_1}&
\dots & a_te_2(x_{2t})x_{2t}^{\alpha_t}
\end{array}\right),
$$
so that $\det(M_0)\equiv \Delta_p\, f(x_1,\dots,x_{2t}) \not\equiv 0
\pmod p$ by the second hypothesis of the theorem.

\smallskip
Since $\det(M_k)\equiv \det(M_0)\not\equiv 0 \pmod p$, $M_k$ is
 invertible and System (\ref{sistema}) has a unique solution
 modulo $p^{2^k}$, namely:
$$\left(\begin{array}{c}d_1\\\vdots\\d_t\\\delta_1\\\vdots\\\delta_t\end{array}
\right)=M_k^{-1}\left(\begin{array}{c}\ell_{2^{k+1}}(x_1)\\\vdots
\\\ell_{2^{k+1}}(x_{2t})\end{array} \right).$$ This shows the existence of
$f_{k+1}$ and its uniqueness property.
\end{proof}

\begin{obs} The statement of Theorem \ref{hensel} also holds modulo $p^k$.
We stated and proved it modulo $p^{2^k}$ to get quadratic
convergence as well as in the classic Newton-Hensel lifting.
\end{obs}

We observe that the hypothesis
$\Delta_p\,f(x_1,\dots,x_{2t})\not\equiv 0 \pmod p$ of Theorem
\ref{hensel} implies in particular that $p\ndiv x_i$, $p\ndiv a_j$
and also that $p-1\ndiv \alpha_j-\alpha_\ell$, $1\le i\le 2t, 1\le
j\ne \ell \le t$ (otherwise, $x_i$ being prime to p would force two
columns of $M_0$ to coincide). Thus $f$ has exactly the same number
$t$ of terms when reduced modulo $p$ than in $\Z[x]$, and no two
exponents reduce to the same modulo $(p-1)$. In view of this we give
the following definition

\begin{defn}\label{reduceswell} Let $p$ be an odd prime number. We say that a polynomial $f=\sum_j
a_jx^{\alpha_j}\in\Z[x]$  with $a_j\ne 0$, $\forall \, j$,  {\em
reduces well modulo} $p$ if $p\ndiv a_j$ for any $j$,  and $p-1\ndiv
\alpha_j-\alpha_\ell$ for any $j\ne \ell$.
\end{defn}

 In these
conditions the uniqueness property of Theorem \ref{hensel} has an
equivalent formulation in terms of how the polynomials coincide as
functions on $(\Z/p^{2^k}\Z)^\times$:

\begin{prop}\label{pkequivalencia}
Let $f=\sum_j a_jx^{\alpha_j},g=\sum_\ell b_\ell
x^{\beta_\ell}\in\Z[x]$ be two polynomials that reduce well modulo
$p$.
 Then, for any $k\in \N$,
 the two following conditions are equivalent:
\begin{itemize}
\item $f$ and $g$ have the same number $t$ of non-zero terms, and
up to an index permutation, $a_j\equiv b_j\pmod{p^k}$ and
$\alpha_j\equiv \beta_j \pmod{\varphi(p^k)}$. \item $f(x)\equiv
g(x) \pmod{p^{k}}$ for all $x\in (\Z/p^k\Z)^\times$.

\end{itemize}
\end{prop}

\smallskip

\smallskip Before giving the proof, we need the following:

\begin{obs}\label{primitive}
Any primitive root modulo $p^{k}$ is also a primitive root modulo $p$.
\end{obs}
\begin{proof} Since $(\Z/p^{k}\Z)^\times$ is a cyclic group,
the number of elements in $(\Z/p^{k}\Z)^\times$ of order divisible by
  $(p-1)$ is $p^{k-1}\varphi(p-1)$. Clearly this set must contain the
  $p^{k-1}$ lifts of any primitive root of $(\Z/p\Z)^\times$ to
  $(\Z/p^{k}\Z)^\times$. By cardinality they are the same.
\end{proof}

\begin{proof}{of Proposition \ref{pkequivalencia}.--}

$(\Downarrow)$ is clear.

$(\Uparrow)$ by induction in $k$: \\
\noindent $\bullet$ Case $k=1$: Write  $f=\sum_{j=0}^{p-2}a_j x^{j
+(p-1)k_j}$ and $g=\sum_{j=0}^{p-2}b_j x^{j+(p-1)k'_j}$. Since
$x^{p-1}\equiv 1 \pmod p$ and $f(x)\equiv g(x) \pmod p$ for $1\le x\le p-1$,
$$\sum_{j=0}^{p-2} (a_j-b_j) x^j\equiv 0 \pmod p\quad \mbox{for}
\ 1\le x\le p-1,$$

but a degree $p-2$ polynomial has at most $p-2$ different roots modulo
$p$ hence $a_j\equiv b_j$ mod$p$.

\smallskip
\noindent $\bullet$ Assume it is true for $k$.

Let $f,g$ be two polynomials satisfying the hypothesis and such that
$f(x)\equiv g(x) \pmod{p^{k+1}}$ for all $x \in
(\Z/{p^{k+1}}\Z)^\times$. In particular $f(x)\equiv g(x) \pmod p$ and
by the $k=1$ case they have both exactly $t$ non-zero terms modulo $p$
for some $1\le t\le p-1$.

If $f:=\sum_{j=1}^ta_jx^{\alpha_j}$ and
$g:=\sum_{j=1}^tb_jx^{\beta_j}$, by inductive hypothesis (up to a permutation of
indexes) $a_j\equiv b_j \pmod{p^k}$ and $\alpha_j\equiv \beta_j \pmod{\varphi(p^k)}$.
Let $a_j^{-1}$  denote the inverse of $a_j$ in the multiplicative
group $(\Z/p^{k+1}\Z)^\times$, and $0\le c_j, \gamma_j<p$ be such
that
$$a_j^{-1}b_j\equiv 1+p^k c_j\pmod{p^{k+1}} \quad \mbox{and} \quad
\beta_j-\alpha_j\equiv \varphi(p^k)\gamma_j
\pmod{\varphi(p^{k+1})}.$$
 For all $x \in (\Z/p^{k+1}\Z)^\times$:
\begin{eqnarray*}
0& \equiv&
(f-g)(x)\pmod{p^{k+1}}\\
&\equiv & \sum_{j=1}^t(a_jx^{\alpha_j}-b_jx^{\beta_j})\pmod{p^{k+1}}\\
&\equiv&\sum_{j=1}^t a_jx^{\alpha_j}(1-a_j^{-1}b_j
x^{{\beta_j}-\alpha_j})\pmod{p^{k+1}}\\
&\equiv&\sum_{j=1}^t a_jx^{\alpha_j}(1-(1+p^kc_j)
x^{\varphi(p^k)\gamma_j})\pmod{p^{k+1}}
\end{eqnarray*}

But $x^{\varphi(p^k)}\equiv 1+ p^ke_{p,k+1,k}(x) \pmod {p^{k+1}}$ for
$e_{p,k+1,k}$ defined in Notation \ref{e-notation}, thus
\begin{eqnarray*}
0&\equiv&\sum_{j=1}^t a_jx^{\alpha_j}(1-(1+p^kc_j)
(1+ p^ke_{p,k+1,k}(x))^{\gamma_j})\pmod{p^{k+1}}\\
&\equiv& \sum_{j=1}^t a_jx^{\alpha_j}(1-(1+p^kc_j)
(1+ p^ke_{p,k+1,k}(x)\gamma_j)\pmod{p^{k+1}}\\
&\equiv& \sum_{j=1}^t -a_jx^{\alpha_j}p^k(c_j+
e_{p,k+1,k}(x)\gamma_j)\pmod{p^{k+1}}
\end{eqnarray*}
Therefore,
\begin{equation}\label{equation1}
 \sum_{j=1}^t a_jx^{\alpha_j}(c_j+
e_{p,k+1,k}(x)\gamma_j)\equiv 0 \pmod{p}.
\end{equation}

Let $\rho\in \Z$ be a primitive root mod$p^{k+1}$. Then $x_1:=1,
x_2:=\rho^{p},\dots, x_{t}:=\rho^{(t-1)p}$ are all different modulo
$p$ by Observation \ref{primitive}, and
$e_{p,k+1,k}(x_i)=e_{p,k+1,k}(\rho^{(i-1)p}) \equiv 0 \pmod p$ since
$e_{p,k+1,k}$ is a group homomorphism. Substituting in
(\ref{equation1}):
$$
 \sum_{j=1}^t a_j\rho^{(i-1)p\,\alpha_j}c_j\equiv 0 \pmod{p}\quad \mbox{for} \ 1\le i\le t.
$$
That is
$$\left( \begin{array}{ccc}
a_1& \dots & a_t \\
a_1\rho^{p\alpha_1}&\dots & a_t\rho^{p\alpha_t}\\
\vdots& & \vdots\\
a_1\rho^{(t-1)p\alpha_1}&\dots &a_t\rho^{(t-1)p\alpha_t}
\end{array}
\right)\left( \begin{array}{c}c_1\\c_2\\ \vdots
\\c_t\end{array}\right) \equiv \left( \begin{array}{c}0\\ \vdots
\\0\end{array}\right) \pmod p.$$

The Vandermonde determinant of this linear system is
 $$a_1\cdots a_t
 \prod_{j<\ell}(\rho^{p\,\alpha_\ell}-\rho^{p\,\alpha_j})\not\equiv 0
 \pmod p,$$ since $\alpha_j\not\equiv \alpha_\ell\mod{(p-1)}$
 . Therefore $c_j\equiv 0 \pmod p$ and $a_j\equiv b_j \pmod {p^{k+1}}$
 for $1\le j\le t$.

 \smallskip
 In a similar way, taking  $x_1:=\rho, x_2:=\rho^2,\dots, x_{t}:=\rho^{t}$
 in (\ref{equation1}) knowing that $c_j=0$ and $e_{p,k+1,k}$ is a group
 homomorphism we get for $e_\rho:=e_{p,k+1,k}(\rho)$:
$$
\sum_{j=1}^t a_j\rho^{i\alpha_j} (i\,e_\rho) \gamma_j\equiv 0
\pmod{p} \quad \mbox{for} \ 1\le i\le t,$$ that leads to a
Vandermonde linear system with  determinant equal to
$$t!\,a_1\cdots a_t\, e_\rho^t \rho^{\alpha_1}\prod_{j<\ell}(\rho^{\alpha_\ell}-\rho^{\alpha_j})
\not\equiv 0 \pmod p,$$ since $e_\rho \not\equiv 0\pmod p$. We
conclude that $\gamma_j\equiv 0 \pmod p$ and therefore
$\alpha_j\equiv \beta_j \pmod {\varphi(p^{k+1})}$.
\end{proof}

The question in order to apply Theorem \ref{hensel} is whether, given
a polynomial $f=\sum_{j=1}^ta_jx^{\alpha_j} \in\Z[x]$ and an odd prime
number $p$ such that $f$ reduces well modulo $p$, there exist
starting interpolating sets $\set{x_1,\dots,x_{2t}} \subset \Z$
satisfying the assumption that
$\Delta_p\,f(x_1,\dots,x_{2t})\not\equiv 0 \pmod p$.

The condition is independent from the polynomial coefficients
$a_i$ since we can factor out $a_1\cdots a_t$ from $\Delta_p\,f$
and $f$ reducing well modulo $p$ implies that $a_i \not \equiv 0
\pmod p$. Proposition \ref{startingset} below gives examples of
good starting sets modulo $p$ independent from  $f$ (i.e. they are
good for any polynomial $f\in\Z[x]$ with exactly $t$ terms that
reduces well modulo $p$).

\begin{defn} \label{defstartingset}Let $p$ be an odd prime number, and
 $t\in \N$, $1\le t\le p-1$.  We say that $\set{x_1,\dots,x_{2t}}
\subset \Z$ is a {\em good starting set modulo $p$} (for
Newton-Hensel interpolation) if for every subset
$\{\alpha_1,\dots, \alpha_t\}\,\subseteq \{0,\dots,p-2\}$, the
quantity \begin{equation}\label{determinant}
 \det \left(
\begin{array}{cccccc}x_1^{\alpha_1} & \dots &
x_1^{\alpha_t}&e_2(x_1)x_1^{\alpha_1}& \dots &
e_2(x_1)x_1^{\alpha_t}\\
 \vdots  & & \vdots  & \vdots & & \vdots \\
x_{2t}^{\alpha_1} &  \dots &
x_{2t}^{\alpha_t}&e_2(x_{2t})x_{2t}^{\alpha_1}&
 \dots & e_2(x_{2t})x_{2t}^{\alpha_t}
\end{array}\right) \not\equiv 0 \pmod p,
\end{equation}
  where $e_2:\Z \rightarrow \Z/p^2\Z \rightarrow \Z/p\Z$
  is the map of Notation \ref{e-notation}.
\end{defn}

\begin{rem} While the previous determinant depends on the order
chosen for the points $x_i$ and the exponents $\alpha_j$, the
condition of being zero (respectively non-zero) does not.
\end{rem}

\begin{prop} \label{startingset}
The following sets $\set{x_1,\dots,x_{2t}}\subset \Z$ are good
starting sets
 modulo $p$:
\begin{enumerate}
\item $\{x_1,\dots,x_{2t}\}$ where $x_i\equiv \rho^{i-1}\pmod{p^2}$, $1\le i\le 2t$, for $\rho\in \Z$
 a primitive root modulo
$p^2$, or, more generally,  where $x_i\equiv \rho^{a+i-1}
\pmod{p^2}$, $1\le i\le 2t$, for $\rho \in \Z$ a primitive root
modulo $p^2$ and $a \in \Z$ ($2t$ consecutive powers of a primitive
root modulo $p^2$).
\item
 $\{x_1,\dots,x_{2t}\}$ where $x_i\equiv\rho^{i-1}\pmod p$ and $x_{t+i}=x_i+p$,
 $1\le i\le t$, for $\rho\in \Z$
 a primitive root modulo $p$, or , more generally,
 where $x_i, x_{t+i}\equiv\rho^{a+i-1}\pmod p$ and $x_{t+i}\not\equiv x_{i}\pmod{p^2}$,
 $1\le i\le t$, for $\rho\in \Z$
 a primitive root modulo $p$ and $ a\in \Z$ (any set
 $x_1,\dots,x_{2t}$ where you choose 2  sets of $t$ consecutive powers
 of a primitive
root modulo $p$, formed by different elements modulo $p^2$).
\end{enumerate}
 \end{prop}

\begin{proof}
We show the main cases of the two items, since their generalizations are straight-forward.

We denote by $\Delta_p$ the determinant modulo $p$ defined in
(\ref{determinant}).

\smallskip
In the first case $x_i\equiv \rho^{i-1}$ modulo $p^2$ for $1\le i\le
2t$. Denote $e:=e_2(\rho)$, then $e \ne 0$ since $\rho$ is a
primitive root modulo $p^2$ and
$e_2(x_i)=e_2(\rho^{i-1})=(i-1)e_2(\rho)=(i-1)e$. Substituting in
the definition we get
\begin{eqnarray*}
\Delta_p =& \det\left( \begin{array}{cccccc}
1&\dots & 1&0& \dots &0\\
\rho^{\alpha_1} &\dots & \rho^{\alpha_t}&e\rho^{\alpha_1} &
\dots &e\rho^{\alpha_t}\\
\rho^{2\alpha_1} & \dots & \rho^{2\alpha_t}&2e\rho^{2\alpha_1} &
\dots &2e\rho^{2\alpha_t}\\
 \vdots &  & \vdots & \vdots &  & \vdots \\
\rho^{(2t-1)\alpha_1} & \dots &
\rho^{(2t-1)\alpha_t}&(2t-1)e\rho^{(2t-1)\alpha_1} & \dots &
(2t-1)e\rho^{(2t-1)\alpha_t}
\end{array}\right)
\\[4mm]
=&  e^t\det\left(
\begin{array}{cccccc}1&
\dots & 1&0&  \dots &
0\\
\rho^{\alpha_1} & \dots & \rho^{\alpha_t}&\rho^{\alpha_1} &
\dots &\rho^{\alpha_t}\\
\rho^{2\alpha_1} & \dots & \rho^{2\alpha_t}&2\rho^{2\alpha_1} &
\dots &2\rho^{2\alpha_t}\\
 \vdots &  & \vdots & \vdots &   & \vdots \\
\rho^{(2t-1)\alpha_1} & \dots &
\rho^{(2t-1)\alpha_t}&(2t-1)\rho^{(2t-1)\alpha_1} & \dots &
(2t-1)\rho^{(2t-1)\alpha_t}
\end{array} \right)
\end{eqnarray*}

If we denote $z_i:=\rho^{\alpha_i}$ for $1\le i \le t$, then

$$ \Delta_p =  e^t\det\left(
\begin{array}{cccccc}1&
\dots & 1&0&\dots &
0\\
z_1 & \dots & z_t&z_1 &
\dots &z_t\\
z_1^2 & \dots & z_t^2&2z_1^2 &
\dots &2z_t^2\\
 \vdots &  & \vdots & \vdots &   & \vdots \\
z_1^{2t-1}& \dots & z_t^{2t-1}&(2t-1)z_1^{2t-1} & \dots &
(2t-1)z_t^{2t-1}
\end{array}  \right)\pmod p$$

The transpose of this last matrix is well-known. It arises in the
 Hermite interpolation problem while trying to interpolate a
function by a polynomial $f$ such that $f(z_i)=y_i$ and
$f'(z_i)={y_i'}$ for $1\le i\le
 t$.  Its determinant equals
$$(-1)^{t(t-1)/2}z_1\cdots z_t \prod_{1\le i<j\le t}(z_j-z_i)^4.$$
Hence
$$\Delta_p = (-1)^{t(t-1)/2}e^t \rho^{\alpha_1+\cdots + \alpha_t}
\prod_{1\le i<j\le t} (\rho^{\alpha_j}-\rho^{\alpha_i})^4 \pmod p.$$
Since $1\le \alpha_i \le p-2$ are all distinct and $\rho$ is a
primitive root modulo $p$ (by Observation \ref{primitive}),
$\rho^{\alpha_i} \not\equiv \rho^{\alpha_j} \pmod p$ if $i \neq j$,
hence $\Delta_p \not\equiv 0 \pmod p$ as wanted.

The general case reduces to this one factoring out from the determinant
the (non-zero) term $\rho^{2a(\alpha_1+\dots +\alpha_t)}$.

\medskip
In the second case $x_i\equiv \rho^{i-1}$ modulo $p$ and
$x_{t+i}=x_i+p$ for $1\le i\le t$. Then $e_2(x_{t+i})\equiv
e_2(x_i)-x_i^{-1} \pmod p$ for $1\le i\le t$ since
$$(x_i+p)^{p-1}\equiv x_i^{p-1} +(p-1)\,p\,x_i^{p-2} \equiv 1+p\,(e_2(x_i)- x_i^{-1}) \pmod{p^2}.$$
Thus, since $x_{t+i}\equiv x_i\pmod p$, we have  modulo $p$:
$$\begin{array}{rcll} \Delta_p &=&\det \left(
\begin{array}{cccccc}x_1^{\alpha_1} & \dots &
x_1^{\alpha_t}&e_2(x_1)x_1^{\alpha_1}& \dots &
e_2(x_1)x_1^{\alpha_t}\\
\vdots & & \vdots &\vdots& &\vdots\\x_t^{\alpha_1} & \dots &
x_t^{\alpha_t}&e_2(x_t)x_t^{\alpha_1}& \dots &
e_2(x_t)x_t^{\alpha_t}\\
x_1^{\alpha_1}&
 \dots & x_1^{\alpha_t}  &(e_2(x_1)-x_1^{-1})x_1^{\alpha_1}& \dots &
 (e_2(x_1)-x_1^{-1})x_1^{\alpha_t} \\
\vdots & &
\vdots & \vdots & & \vdots\\
x_t^{\alpha_1}&
 \dots & x_t^{\alpha_t}  &(e_2(x_t)-x_t^{-1})x_t^{\alpha_1}& \dots &
 (e_2(x_t)-x_t^{-1})x_t^{\alpha_t}
\end{array}\right) &  \pmod{p}
\\[1.5cm]
&= & \det\left(  \begin{array}{cccccc}x_1^{\alpha_1} & \dots &
x_1^{\alpha_t}&e_2(x_1)x_1^{\alpha_1}& \dots &
e_2(x_1)x_1^{\alpha_t}\\
\vdots & & \vdots &\vdots&
&\vdots\\
x_t^{\alpha_1} & \dots & x_t^{\alpha_t}&e_2(x_t)x_t^{\alpha_1}&
\dots &
e_2(x_t)x_t^{\alpha_t}\\
0&
 \dots & 0  &-x_1^{\alpha_1-1}& \dots &
 -x_1^{\alpha_t-1}\\
\vdots & &
\vdots & \vdots & & \vdots\\
0&
 \dots & 0  &-x_t^{\alpha_1-1}& \dots &
 -x_t^{\alpha_t-1}
\end{array}\right) & \pmod{p}\\[1.5cm]
&=& (-1)^tx_1^{-1}\dots x_t^{-1}\prod_{1\le i<j\le
t}(\rho^{\alpha_j}-\rho^{\alpha_i})^2 \ \not\equiv \ 0 & \pmod{p}.
\end{array}$$

\smallskip
For the general case observe that for $1\le i\le t$,
 $x_{t+i}\equiv x_i+k_ip\pmod{p^2}$
for some $k_i$ prime to $p$.
\end{proof}

\smallskip
Unfortunately, not every set $\{x_1,\dots,x_{2t}\}$ with
$x_i\not\equiv 0\pmod{p}$ is a good starting set modulo $p$:

\begin{exmpl}  Take $p:=7$, $\rho:=3$ and $t=2$, $ x_1\equiv 3^0
  \pmod 7$, $x_2\equiv 3^1 \pmod 7$, $x_2:\equiv 3^2 \pmod 7$ and
$x_3\equiv 3^4 \pmod 7$, then $e_2(3)=6\equiv -1 \pmod 7$, for
$\alpha_1:=0$ and $\alpha_2:=3$ we obtain:

$$
\Delta_p\equiv\det \left(
\begin{array}{cccc}
1&1&0 &0\\
3^{0} &  3^{3}&-3^{0} &
-3^{3}\\
3^{2\cdot 0} & 3^{2\cdot 3}&-2\cdot 3^{2\cdot 0} &
-2\cdot 3^{2\cdot 3}\\
3^{4\cdot 0} &
 3^{4\cdot 3}&-4\cdot 3^{4\cdot 0}
& -4\cdot 3^{4\cdot 3}
\end{array}\right)\equiv \det
 \left(
\begin{array}{rrrr}1&1&0 &0\\
1 &  -1&-1 &
1\\
1 & 1&-2 &
-2\\
1 &
 1&-4
& -4
\end{array} \right) \equiv 0 \!\!\pmod 7.
$$
\end{exmpl}

\smallskip
In particular, Proposition \ref{startingset} implies that for any
polynomial $f\in \Z[x]$ that reduces well modulo $p$ and with
exactly $t$ non-zero terms , there exists a subset of size $2t$ in
$\{1,\dots, 2p-1\}$, $x_1,\dots, x_{2t}$, that satisfies the
hypothesis of Theorem \ref{hensel}, i.e. such that
$\Delta_p\,f(x_1,\dots,x_{2t})\not\equiv 0 \pmod p$. In the case of
a binomial, this result can be sharpened: for any prime $p\ge 5$,
there exists a subset $\{x_1,\dots,x_4\}\subset \{1,\dots,p-1\}$
with $\Delta_p f(x_1,\dots,x_4)\not \equiv 0 \pmod p$.

\begin{prop}\label{binomio} Let $p\ge 5$ be a prime number and
$0\le \alpha < \beta <p-1$. Then
$$\left( \begin{array}{cccc}1^{\alpha} &
1^{\beta}&e_2(1)1^{\alpha}&
e_2(1)1^{\beta}\\
2^{\alpha} & 2^{\beta}&e_2(2)2^{\alpha}&
e_2(2)2^{\beta}\\
\vdots & \vdots &\vdots  & \vdots \\ (p-1)^{\alpha} &
(p-1)^{\beta}&e_2(p-1)\,(p-1)^{\alpha}&
 e_2(p-1)\,(p-1)^{\beta}
\end{array}\right) \pmod p
$$
has rank four, where $e_2:\Z \rightarrow \Z/p^2\Z \rightarrow
\Z/p\Z$
  is the map of Notation \ref{e-notation}.
\end{prop}

We omit the proof of this fact since it is quite tedious and is
based on a smart choice of the four elements $x_i$ depending on
congruences of $\alpha - \beta$ modulo $2$, modulo $(p-1)$ and on
the behavior of $e_2(x_i)$. \hfill\mbox{$\Box$}

\medskip

Unfortunately this result is not true in general since for $p:=11$,
the polynomial $f=a_1+a_2\,x+a_3\,x^3+a_4\,x^5+a_5x^8$ has five
non-zero terms but $\Delta_p\,f(1,\dots,10)\equiv 0\pmod{11}$.
Surprisingly, computer experiments performed with \cite{GP} did not
show any counter-example for trinomials (with any prime) nor
polynomials with $5$ terms for a prime different from $11$.

\smallskip For an arbitrary odd prime number $p$, there may be more
good starting sets modulo $p$ for a given $t$ than the sets
described in Proposition \ref{startingset} but we are still unable
to prove their existence in general. For instance for $p=7$ and
$t=2$, the set $\{1,2,3,6\}$ is good but is neither of type (1) nor
of type (2). Actually in that case the total number of sets is
$10626$, from which  $1640$ are good but only $560$ are of type (1)
or (2).

\smallskip
In what follows, we compute the number of different good starting
sets as in Proposition \ref{startingset} on the interval
$\{1,\dots,p^2-1\}$ in order to estimate (at least from below) the
probability of having a good starting set modulo $p$ when choosing
any subset of $\{1,\dots,p^2-1\}$.

A first question we have to deal with is whether  the sets described
of type (1) and type (2) of Proposition \ref{startingset} are
distinct when choosing arbitrarily a primitive root modulo $p^2$ or
modulo $p$ and a starting exponent $a$.
 To make this analysis we need to distinguish the case $t=p-1$ from $t<p-1$.

\begin{obs} If $t=p-1$, the good starting sets modulo $p$ of type (1)
in Proposition \ref{startingset} are contained in type (2) while if
$t<p-1$ type (1) and (2) define different good starting sets.
\end{obs}

\begin{proof} For $t=p-1$, this is clear since the set of
$p-1$ consecutive powers of a primitive root modulo $p$ coincides with
  the set $\{1,\dots, p-1\}$ and
 a primitive root modulo $p^2$ is a primitive root
modulo $p$. For $t<p-1$, in a set $\{\rho^a, \rho^{a+2t-1}\}$ of
type (1) there are more than $t$ consecutive primitive roots modulo
$p$ that cannot appear in a set of type (2).
\end{proof}

\begin{prop}
For $t=p-1$ there are  at least ${p\choose 2}^{p-1}$  good starting
sets modulo $p$ in the set $\{1,\dots, p^{2}-1\}$.
\end{prop}

\begin{proof} We simply choose for each $i$, $1\le i\le p-1$, two
different  elements congruent to $i$ modulo $p$. There are
${p\choose 2}^{p-1}$ such choices.
\end{proof}

Taking into
 account that the number of subsets of size $2(p-1)$ in $\{1,\dots,p^2-1\}$, avoiding the
 multiples of $p$,  equals ${p(p-1) \choose 2(p-1)}$, this gives a ratio of
 $$\frac{{p\choose 2}^{p-1}}{{p(p-1)\choose 2(p-1)}}\approx \frac{1}{p(p-1)^{p-1}}$$
 good starting sets modulo $p$, which is obviously a very small quantity.

\bigskip
For $t<p-1$ we need to perform a more careful analysis:

\begin{lem} \label{Znmaps}
Let $n\in \N, \, n>2$ and $\CC:=\{0,1,\dots,s\} \subset \Z/n\Z$ where
  $1 \le s <n-1$. If $\sigma(x)=m\,x+a$ is a bijective map of $\Z/n\Z$
  into itself such that $\sigma(\CC) = \CC$, then $\sigma = Id$ or
  $\sigma=- Id + s$.
\end{lem}

\begin{proof}
Without loss of generality we may suppose $1\le m<n$ and $0\le a
<n$. Furthermore we may restrict to the case $s<n/2$ since if not
$\overline{\CC}=\{s+1, \ldots , n-1\}$ is also fixed by $\sigma$. Then
the bijective map $\tau(x) = \sigma(x+s+1)-(s+1)$ fixes $\{0,1,
\ldots, n-s-2\}$ (where $n-s-2<n/2$) and satisfies that $\sigma=Id
\Leftrightarrow \tau = Id$ and $\sigma=-Id +s \Leftrightarrow
\tau=-Id+(n-s-2)$.
We consider two different cases:
\begin{itemize}\item Case $a=0$. Note that in this case $m<n/2$ since
 $\sigma(1)=m\in \CC$.  We need to show that $\sigma=Id$.  Assume that
  $m>1$, then by the Euclidean algorithm $s=qm+r$ with $0\le r<m$. In
  particular $q<s$ since $m>1$, therefore $q+1\in \CC$ and
  $\sigma(q+1)\in \CC$, i.e. $\sigma(q+1)\le s$.  But $0\le (q+1)m\le
  s+m<n/2+n/2=n$ implies that $\sigma(q+1)=(q+1)m$, and, on the other
  side, $(q+1)m>qm+r=s$. Contradiction.
\item Case $a>0$. Assuming $m >1$ we will prove that $m=n-1$ and $a=s$.

Note that $m > s$ since if not let $\AA:=\{0\le x\le n-1: mx+a\le
s\}$. Clearly $\AA$ is non-empty ($0\in \AA$), denote by $y$ its
maximum element, i.e. $m y+a\le s$ and $m(y+1)+a>s$. Since $m \le s$,
$0\le m(y+1)+a\le 2s<n$, then $\sigma(y+1)=m(y+1)+a$. Since $m>1$ and
$m y+a\le s$, $y<s$ thus $y+1\in \CC$. But then $\sigma(y+1)\le s$ and
also $m(y+1)+a >s$, a contradiction. Then $s<m\le n-1$.

In that case, $\sigma(s)=0$ since if $\sigma(y)=0$ for $y<s$, then
$\sigma(y+1)=m\le s$ ($\sigma $ is bijective on $\CC$), which is
not the case. Thus there exists $z<s$ such that $\sigma(z)=1$ and
$\sigma(z+1) = m+1 \in \CC$. That is, $s+1<m+1\le n $ belongs to $\CC$
mod$n$. That means that $m+1\equiv 0$ mod $n$ and $m=n-1$. We
conclude, since $0=\sigma(s)=-s+a$, that $a=s$.
\end{itemize}
\end{proof}

\begin{cor} Let $\rho, \varrho \in \Z$ be two different primitive roots mod$p^k$.
  Then, for any $s$, $1\le
 s<\varphi(p^k)$, the sets
$\rho^a \{1,\rho,
  \ldots, \rho^{s-1}\}$ and $ \varrho^{b} \{1,\varrho, \ldots, \varrho^{s-1}\}$
  coincide modulo $p^k$ if and only if $\varrho \equiv \rho\pmod{p^k}$ and $b \equiv
  a\pmod{\varphi(p^k)}$ or $\varrho \equiv
  \rho^{-1}\pmod{p^k}$ and $b\equiv-(a+s-1)\pmod{\varphi(p^k)}$.
\end{cor}

\begin{proof} Without loss of generality we may assume
  $\varrho = \rho^m$ for some  $m$ prime to $\varphi(p^k)$. Then  the
  two sets coincide if and only if $\{1,\rho,
  \ldots, \rho^{s-1}\}$ and $ \rho^{mb-a} \{1,\rho^m, \ldots, \rho^{m\,{s-1}}\}$
  coincide, i.e. if the bijective
  map $\sigma(x) = m\,x+(mb-a)$ of $\Z/\varphi(p^k)\Z$ fixes the set $\{0,1, \ldots,s-1\}$.
   We apply  Lemma \ref{Znmaps} to conclude that
   $\sigma=Id$, that is $m\equiv 1\pmod{\varphi(p^k)}$,  and $a\equiv
   b\pmod{\varphi(p^k)}$ or $\sigma=-Id+s$, that is
   $m\equiv -1\pmod{\varphi(p^k)}$ and $b\equiv -(a+s-1)\pmod{\varphi(p^k)}$.
\end{proof}

As an immediate consequence we can compute the number of good sets
of type (1) and (2) in Proposition \ref{startingset} for $t<p-1$:

\begin{cor} Let $p$ be an odd prime number and set $1\le t< p-1$.
There are at least $$p\,\varphi(p-1)(p-1)^2/2 + {p\choose
2}^t\varphi(p-1)(p-1)/2$$   good starting sets modulo $p$ in the set
$\{1,\dots, p^2-1\}$.
\end{cor}

\begin{proof}
First we count  the sets of
 type (1) in Proposition \ref{startingset}:
We choose a primitive root $\rho$ modulo $p^2$, there are $\varphi(\varphi(p^2))$ such choices,
and  we choose $0\le a<\varphi(p^2)$. Since $t<p-1$, $2t<\varphi(p^2)$  we can apply the previous lemma
and divide the total number by 2.

Next we count the sets of type (2): We first compute the number of
sets modulo $p$, i.e. the number of different sets  $\set{
\rho^{a+i-1} \pmod p \st 1\le i\le t}$
 in   the set $\{1,\dots,p-1\}$. Applying the previous lemma,
 there are $\varphi(p)\varphi(\varphi(p))/2$ such different choices.
 Then for $1\le i\le t$ we freely choose $x_i,x_{t+i}$ congruent to $\rho^{a+i-1}$
 but different. There are ${p\choose 2}^t$ such choices.
\end{proof}

 Like in the case $t=p-1$, the ratio is again very deceptive: given $p$
 and $t$ the probability of
 choosing randomly a good starting set modulo $p$ is very low and tends to $zero$
 when $p$ grows. However, a more realistic probability estimation
 would be to  compute, given an odd prime number $p$ and $f\in \Z[x]$ a
 polynomial with exactly $t<p$ non zero terms that reduces well
 modulo $p$,  the probability that randomly chosen $\{x_1,\dots,
 x_{2t}\}\subset\{0,\dots, p^2-1\}$ satisfy that Determinant
 (\ref{determinant}) does not vanish modulo $p$. Unfortunately we are still not able to
 give a sharp estimation for that probability.

\section{$p$-adic equations}
\label{$p$-adic equations}

This section shows how  Newton-Hensel construction of Theorem
\ref{hensel} corresponds to an usual Hensel lemma on the $p$-adic
integers $\Z_p$. It  explains the role played by the $e_2$ group
homomorphism of Proposition \ref{e-morphism} and why  starting sets
modulo $p$ need to contain a primitive root modulo $p^2$. We begin
by recalling some definitions and properties of $\Z_p$. We refer to
\cite[Ch.II]{Serre70} for the details.

\smallskip
For a prime integer $p$, the set of  $p$-adic integers
$\Z_p$ is the inverse limit of the diagram
$$ \Z/p\Z \,\mathop{\longleftarrow}^{\phi_{1}}\, Z/p^{2}\Z
\,\mathop{\longleftarrow}^{\phi_{2}}\,
 \Z/p^{3}\Z\mathop{\longleftarrow}^{\phi_{3}}\, \cdots \, ,$$
where $\phi_k:\Z/p^{k+1}\Z \to \Z/p^k\Z$ is the canonical projection.
 Here we will view   $\Z_p$ as the equivalent construction
$$\mathbb{Z}_{p}=\left\{ (a_{k})_{k\in \N}\in\mathbb{Z}^{\mathbb{N}}\; :\; a_{k+1}\equiv
a_{k}\pmod{p^{k}}\ \forall\, k\in \N\,\right\} /\sim $$ where
$\sim$ is the equivalence relation defined by
 $$(a_{k})_{k\in \N}\sim(b_{k})_{k\in \N}\ \Longleftrightarrow \ a_{k}\equiv
b_{k}\pmod{p^{k}}\;\forall\, k\in \N.$$ and the  operations are
coordinate-wise, i.e.:
$$\overline{(a_{k})}_{k\in\N}+\overline{(b_{k})}_{k\in \N}:=
\overline{(a_{k}+b_{k})}_{k\in \N}\quad \mbox{and} \quad \overline{(a_{k})}_{k\in \N}\cdot
\overline{(b_{k})}_{k\in \N}:=\overline{(a_{k}\cdot b_{k})}_{k\in \N}.$$
In this formulation,  $a:=\overline{(a_k)}_{k\in \N} \in \Z_p^\times$,
 the multiplicative group of $p$-adic units, if and only if
$ p \ndiv a_1$.

\smallskip From now on $p$ will denote an odd prime number.  In this
case the multiplicative group $\Z_p^\times$ and the additive group
$\Z/(p-1)\Z\,\times \,\Z_p$ are isomorphic. This last additive
group can be viewed in another additive way, closely related to
the changes of exponents we allow in Theorem \ref{hensel} and that
justifies a convenient exponentiation in $\Z_p^\times$. For that
purpose, we construct the ring of {\em exponents} that we will
denote $\EE_p$:

\begin{defn} \label{Ep} Let $p$ be an odd prime number. We define $\EE_p$ as the inverse limit of the
following diagram:
$$ \Z/\varphi(p)\Z \,\mathop{\longleftarrow}^{\rho_{1}}\, Z/\varphi(p^{2})\Z
\,\mathop{\longleftarrow}^{\rho_{2}}\,
 \Z/\varphi(p^{3})\Z \mathop{\longleftarrow}^{\rho_{3}}\, \cdots \, ,$$
where, since $\varphi(p^{k})\div \varphi(p^{k+1})$,
 $\rho_k:\Z/\varphi(p^{k+1})\Z \to \Z/\varphi(p^k)\Z$ is the well-defined
  canonical projection.
\end{defn}
Equivalently, $\EE_p$ can be seen as the ring
$$\EE_{p}=\left\{ (\alpha_{k})_{k\in \N}\in\Z^{\N}\; :\; \alpha_{k+1}\equiv
\alpha_{k}\pmod{\varphi(p^{k})}\ \forall\, k\in \N\right\}
/\approx$$ where
 $\approx$ is the equivalence relation defined by
 $$(\alpha_{k})_{k\in \N}\approx (\beta_k)_{k\in \N}\ \Longleftrightarrow \
\alpha_k\equiv \beta_k\pmod{\varphi(p^{k})}\;\forall\, k\in \N.$$
where the  operations are coordinate-wise, i.e.:
$$\overline{(\alpha_k)}_{k\in \N}+\overline{(\beta_k)}_{k\in \N}:=
\overline{(\alpha_k+\beta_{k})}_{k\in \N}\quad \mbox{and} \quad \overline{(\alpha_k)}_{k\in \N}\cdot
\overline{(\beta_k)}_{k\in \N}:=\overline{(\alpha_k\cdot \beta_k)}_{k\in \N}.$$

In this formulation, $\alpha:=\overline{(\alpha_k)}_{k\in \N}\in \EE_p^\times$ if and only if
$\gcd(\alpha_2,\varphi(p^2))=1$, and we have the ring isomorphism
 $$\EE_p \simeq \Z/(p-1)\Z \times \Z_p: \quad \overline{(\alpha_k)}_{k\in \N}
\mapsto (\alpha_1, \overline{(\alpha_k)}_{k\ge 2}).  $$

\smallskip
Next proposition shows that $\EE_p$ is a natural set of exponents for $\Z/p\Z^\times$:

\begin{prop} Let $a:=\overline{(a_k)}_{k\in \N} \in \Z_p^\times$ and
$\alpha:=\overline{(\alpha_k)}_{k\in \N} \in \EE_p$, then $a^\alpha:=
\overline{(a_k^{\alpha_k})}_{k\in \N}\in \Z_p^\times$.
\end{prop}

\begin{proof}
{\ }
\begin{itemize}
\item It is immediate that $a_k\equiv a'_k\pmod{p^k}$ and
 $\alpha_k\equiv \alpha'_k\pmod{\varphi(p^k)}$
implies $a_k^{\alpha_k}\equiv (a'_k)^{\alpha_k}\equiv
(a'_k)^{{\alpha'}_k}\pmod{p^k}$ since $a_k\not\equiv 0\pmod{p}$.
\item Similarly, $a_{k+1}\equiv a_k\pmod{p^k}$ and
 $\alpha_{k+1}\equiv \alpha_k\pmod{\varphi(p^k)}$
implies $a_{k+1}^{\alpha_{k+1}}\equiv a_{k+1}^{\alpha_{k}}\equiv
a_k^{\alpha_k}\pmod{p^k}$.
\end{itemize}
\end{proof}

\begin{cor} For a given  $\alpha\in \EE_p$, the map
$\Z_p^\times \to \Z_p^\times: a\mapsto a^\alpha$ is a group
homomorphism. Moreover if $\alpha\in \EE_p$ is invertible, this
map is an isomorphism since $b=a^\alpha \Leftrightarrow
a=b^{\alpha^{-1}}$.
\end{cor}

\begin{lem} \label{goodbasis} For a given  $a\in \Z_p^\times$, the map
$(\EE_p ,+)\to (\Z_p^\times,\cdot): \alpha \mapsto a^\alpha$ is a
group homomorphism. Moreover, this map is an isomorphism if and
only if $a:=\overline{(a_k)}_{k\in \N}$ is such that $a_2$ is a
primitive root modulo $p^2$. We call such $a$ a {\em good basis
for taking logarithms}.
\end{lem}

\begin{proof} Note that if $a:=\overline{(a_k)}_{k\in \N}$
is such that $a_2$ is a primitive root modulo $p^2$, then $a_k$ is a
primitive root modulo $p^k$ for all $k\in \N$ (see for instance
\cite[Th.10.6]{Apostol76}). Now for $b:=\overline{(b_k)}_{k\in
\N}\in \Z_p^\times$, we let $\alpha_k $ be such that
$a_k^{\alpha_k}\equiv b_k\pmod{p^k}$ and we define
$\alpha:=\overline{(\alpha_k)}_{k\in \N}$. It is easy to check that
$\alpha\in \EE_p$, and  $a^\alpha=b$.

Conversely, we want to check that $a_2$ is a primitive root modulo
$p^2$ if the map is  onto. For each $1 \le b_2 < p^2$, with
$\gcd(b_2,p)=1$, let $b:=\overline{(b_2)}_{k\in \N}\in \Z_p^\times$
(the natural injection of $\Z$ into $\Z_p$),  and let $\alpha$ be
such that $a^\alpha=b$. In particular $a_2^{\alpha_2}\equiv
b_2\pmod{p^2}$, i.e. the powers of $a_2$ span $(\Z/p^2\Z)^\times$
then $a_2$ is a primitive root.
\end{proof}

In particular,  if $\rho\in \Z$ is a primitive root modulo $p^2$, then
$\rho:=\overline{(\rho)}_{k\in \N}\in \Z_p^\times$ is  a
good basis for taking logarithms.

\bigskip In view of the previous discussion, for $a_j\in
\Z_p^\times, \alpha_j\in \EE_p, 1\le j\le t$, the specialization
$$\Z_p^\times \rightarrow \Z_p: \quad
 x \mapsto f(x):=\sum_{j=1}^t a_j
x^{\alpha_j} $$ is a well-defined map. As a consequence of
Proposition \ref{pkequivalencia}, two such maps $f=\sum_{j=1}^t
a_jx^{\alpha_j}$ and $g= \sum_{j=1}^t b_jx^{\beta_j}$ coincide on
$\Z_p^\times$ if and only if $a_j=b_j$ in $\Z_p^\times$ and
$\alpha_j=\beta_j$ in $\EE_p$, $1\le j\le t$, i.e. $f=g$. In this
setting Theorem \ref{hensel} admits the following formulation:

\begin{teo} \label{zphensel}{\em (Newton-Hensel interpolation lifting on $\Z_p$)}

Let $p$ be an odd prime number, and $t\in \N$, $1\le t\le {p-1}$.
Set $y_1,\dots, y_{2t}\in \Z_p$. Let
$f_1:=\sum_{j=1}^ta_{j1}x^{\alpha_{j1}}\in\Z[x]$ and
$x_1=\overline{(x_{1k})}_k,\dots,x_{2t}=\overline{(x_{(2t)k})}_k\in
\Z_p^\times$ satisfy

\begin{itemize}
\item $f_1(x_i)\equiv y_i \pmod{p}$, $1\le i\le 2t$,
\item $\det
\left( \begin{array}{ccclcl}x_1^{\alpha_1} & \dots &
x_1^{\alpha_t}&a_1e_2(x_{12})x_1^{\alpha_1}& \dots &
a_te_2(x_{12})x_1^{\alpha_t}\\
 \vdots  & & \vdots  & \vdots & & \vdots \\
x_{2t}^{\alpha_1} & \dots &
x_{2t}^{\alpha_t}&a_1e_2(x_{(2t)2})x_{2t}^{\alpha_1}& \dots &
a_te_2(x_{(2t)2})x_{2t}^{\alpha_t}
\end{array}\right)\not\equiv 0\pmod p$,

where $e_2(x_i):(\Z/p^2\Z)^\times \rightarrow \Z/p\Z$ is the group
homomorphism of Definition \ref{morphism}.

\end{itemize}
Then there exists a
 unique $f:= \sum_{j=1}^t a_jx^{\alpha_j}$ with $a_j\in \Z_p^\times $ and $\alpha_j\in \EE_p$,
  $1\le j\le t$,
that satisfies simultaneously:
\begin{itemize}
\item $f(x_i)=y_i$ in $\Z_p$,  $1\le i\le 2t$.
\item $
a_j=\overline{(a_{j1},\dots )},
\alpha_j=\overline{(\alpha_{j1},\dots )}$, $1\le j \le t$.
\end{itemize}
\end{teo}

\bigskip
In fact, this theorem is a particular case of a more general type
of equations in $\Z_p$. We deal now with two types of equations:
 generalized polynomial ones
 and exponential ones, and we exhibit the duality they
inherit from Lemma \ref{goodbasis}.

\begin{defn}

{\ }
\begin{itemize}
\item A {\em generalized polynomial equation in $\Z_p^\times$}
 is $$f(x):=\sum_{j=1}^ta_jx^{\alpha_j}$$ where
  $t\in \N ,a_j\in \Z_p$ and $\alpha_j\in \EE_p$,  $1\le j\le t$,  are given, and
  $x\in \Z_p^\times$ is the unknown.
\item An {\em exponential equation in $\EE_p$}
 is $$g(\alpha):=\sum_{j=1}^ta_jx_j^{\alpha}$$ where
  $t\in \N, a_j\in \Z_p$  and $x_j\in \Z_p^\times$, $1\le j\le t$, are given, and
  $\alpha\in \EE_p$ is the unknown.
  \end{itemize}
\end{defn}

\begin{obs} \label{dual} Solving a  generalized polynomial equation in $\Z_p^\times$ or an exponential equation in
$\EE_p$ are essentially the same problem, since if $a\in \Z_p^\times$ is a good basis for
taking logarithms, setting  $x=a^\alpha$ and $x_i=a^{\alpha_i}$ we obtain
$$\sum_{j=1}^ta_j\,x^{\alpha_j}= \sum_{j=1}^ta_j\,x_j^{\alpha}.$$
\end{obs}

The two next propositions generalize Hensel's lemma, traditionally  for polynomials in $\Z[x]$ or
$\Z_p[x]$,
to the previous equations. We state them before giving their proofs to highlight their
dual character.

\begin{prop} \label{henselpolyequation}
 Set  $t\in \N, \,y,\, a_j\in \Z_p$
and $\alpha_j:=\overline{(\alpha_{jk})}_k\in \EE_p$, $1\le j\le t$.
Let $x_1\in  \Z$, $p\ndiv x_1$,  be such that
\begin{itemize}
\item $f(x_1):=\sum_{j=1}^ta_j\,x_1^{\alpha_j}\equiv y \pmod p$,
\item $\sum_{j=1}^t a_{j}\,\alpha_{j2}\,x_1^{\alpha_{j}-1}\not\equiv 0 \pmod p$,
\end{itemize}
then there exists a unique $x=\overline{(x_1,x_2,\dots )}\in \Z_p^\times$ such that
$f(x)=y$ in $\Z_p$.
\end{prop}

\begin{prop} \label{henselexpequation}
 Set  $t\in \N,\, y, \,a_j\in \Z_p$
and $x_j:=\overline{(x_{jk})}_k\in \Z_p^\times$, $1\le j\le t$.
Let $\alpha_1\in \Z$ be such that
\begin{itemize}
\item $g(\alpha_1):=\sum_{j=1}^ta_j\,x_j^{\alpha_1}\equiv y \pmod
p$, \item $\sum_{j=1}^t
a_{j}\,e_2(x_{j2})\,x_{j}^{\alpha_{1}}\not\equiv 0 \pmod p$ where
$e_2$ is the group homomorphism defined in
Definition~\ref{morphism},
\end{itemize}
then there exists a unique $\alpha=\overline{(\alpha_1,\alpha_2,\dots)}\in \EE_p$ such that
 $g(\alpha)=y$ in $\Z_p$.
\end{prop}

\begin{proof}{(of Proposition \ref{henselpolyequation})}

We are looking for $x:=\overline{({x_k})}_k \in \Z_p^\times$ such that
$$f(x_k)=\sum_{j=1}^ta_j\,x_k^{\alpha_j}\equiv \sum_{j=1}^t
a_{jk}\,x_k^{\alpha_{jk}} \equiv y_k \pmod{p^k}.$$ We construct it
inductively starting from the given $x_1$. The condition $p\ndiv
x_1$ guarantees that $x\in \Z_p^\times$. Let $k\in \N$ and assume
there is a unique $x_k$ such that $x_k\equiv x_1\pmod{p}$ and
$f(x_k)\equiv y_k \pmod {p^k}$. If we denote
$a_j:=\overline{(a_{jk})}_k$, $y:=\overline{(y_k)}_k$, we are
looking for $x_{k+1}$ such that $f(x_{k+1})\equiv y_{k+1} \pmod
{p^{k+1}}$. This implies $f(x_{k+1})\equiv y_k\pmod{p^k}$ and thus
$x_{k+1}= x_k+p^k\xi$ where $\xi$ is to be uniquely determined
modulo $p$. We use the same arguments as in the proof of Theorem
\ref{hensel}. Since $a_{j(k+1)}\equiv a_{jk} \pmod {p^k}$,
$y_{k+1}\equiv y_k \pmod {p^k}$ and $\alpha_{j(k+1)}\equiv
\alpha_{jk} \pmod {\varphi(p^k)}$, using the fact that $f(x_k)\equiv
y_k \pmod {p^k}$  there exists $z\in \Z$ such that

$$ y_{k+1}-\sum_{j=1}^t a_{j(k+1)} x_k ^{\alpha_{j(k+1)}} \equiv p^k \,z \pmod{p^{k+1}}.$$

We obtain, by Newton expansion,
$$\begin{array}{crclc}&f(x_{k+1})&\equiv & \sum_{j=1}^t a_{j(k+1)}(x_k+p^k\xi)^{\alpha_{j(k+1)}}
& \pmod{p^{k+1}}\\[1mm]&&\equiv&\sum_{j=1}^t ( a_{j(k+1)}
x_k^{\alpha_{j(k+1)}}+p^k\,\alpha_{j(k+1)}\,x_k^{\alpha_{j(k+1)-1}}\,\xi)
&\pmod{p^{k+1}}\\[1mm]
&&\equiv&  y_{k+1} &\pmod{p^{k+1}}\\[1mm]
\iff&p^k\,z &\equiv & p^k\,\sum_{j=1}^t a_{j(k+1)}
\,\alpha_{j(k+1)}\,x_k^{\alpha_{j(k+1)}-1}\xi &\pmod{p^{k+1}}\\[1mm]
\iff&z &\equiv & (\sum_{j=1}^t a_{j(k+1)}
\,\alpha_{j(k+1)}\,x_k^{\alpha_{j(k+1)}-1})\,\xi &\pmod p.
\end{array}
 $$

Now, by hypothesis, since $\alpha_{j(k+1)}\equiv
\alpha_{j2}\pmod{p}$ for $p\div \varphi(p^2)$, we conclude:

$$\begin{array}{rcll}0&\not\equiv& \sum_{j=1}^t a_{j}
\,
\alpha_{j2}\,x_1^{\alpha_{j}-1}& \pmod p\\[2mm]
0&\not\equiv & \sum_{j=1}^t a_{j(k+1)} \,
\alpha_{j(k+1)}x_1^{\alpha_{j(k+1)}-1} &\pmod{p},
\end{array}
$$

and therefore, there exists a unique $\xi$ modulo $p$ that solves the
problem.
\end{proof}

\begin{proof}(of Proposition \ref{henselexpequation})
The exponential equation $g(\alpha)=\sum_{j=1}^ta_jx_j^\alpha=y$ has
a unique solution $\alpha\in \EE_p$ such that $\alpha\equiv
\alpha_1\pmod{(p-1)}$ if and only if the generalized polynomial
equation $f(\xi)=\sum_{j=1}^ta_j\xi^{\beta_j}=y$ has a unique
solution $\xi\in \Z_p^\times$ such that $\xi\equiv \xi_1\pmod{p}$,
where $\xi:=b^\alpha$ and $x_j:=b^{\beta_j}$, $1\le j\le t$, for a
good basis $b$ for taking logarithms.

We check the assumption of Proposition \ref{henselpolyequation}:
   \begin{eqnarray*}
\sum_{j=1}^t a_{j1}\, \beta_{j2}\,\xi_1^{\beta_{j1}-1}\not\equiv 0
\pmod p &\Longleftrightarrow&
\sum_{j=1}^t a_{j1} \,\beta_{j2}\,b_1^{\alpha_1(\beta_{j1}-1)}\not\equiv 0\pmod p\\
&\Longleftrightarrow&\sum_{j=1}^t a_{j1}
\,\beta_{j2}\,x_{j1}^{\alpha_1}\not\equiv 0\pmod p.
\end{eqnarray*}
Here we used  $b_1^{-\alpha_1}\not\equiv 0\pmod{p}$ since
 $p-1\ndiv \alpha_1$ and $b_1$ is a primitive
root modulo $p$.

Now, we apply the group homomorphism $e_2$ to
$x_{j2}=b_2^{\beta_{j2}}$:
 $e_2(x_{j2})\equiv\beta_{j2}\,e_2(b_2)\pmod{p}$. Moreover, since $b_2$ is a primitive root
 mod$p^2$, $e_2(b_2)\not\equiv 0\pmod{p}$.
 Thus
$$\sum_{j=1}^t a_{j1}\, \beta_{j2}\,x_{j1}^{\alpha_1}\not\equiv 0\pmod p
\Longleftrightarrow \sum_{j=1}^t
a_{j1}\,e_2(x_{j2})\,x_{j1}^{\alpha_1}\not\equiv 0\pmod p.$$
\end{proof}

We observe that in both propositions, the second assumption on $x_1$
and $\alpha_1$ respectively is the one that corresponds to the usual
$f'(x_1)\not\equiv 0 \pmod p$ in the classical Newton-Hensel
lemma. Here, taking derivatives is not a closed operation (because
$\alpha_2\not\equiv \alpha_1 \pmod p$ but modulo $(p-1)$) and the
assumption on the second order is the natural replacement of the
derivative.  Therefore, both for a generalized polynomial equation
$f(x)=\sum_ja_jx^{\alpha_j}$ and for an exponential equation
$g(\alpha)=\sum_ja_jx_j^{\alpha}$, we call the expressions
  \begin{equation}\label{pseudo-derivative}
  \Delta_p\,f(
 x):=\sum_ja_{j}\,\alpha_{j2}\,x^{\alpha_{j}-1}\pmod p\quad \mbox{and} \quad
\Delta_p\,g(
 \alpha):=\sum_ja_{j}\,e_2(x_{j2})\,x_{j}^{\alpha_{1}} \pmod
 p
 \end{equation}
 their  {\em pseudo-derivatives} modulo $p$.

 \bigskip
Like the usual Newton-Hensel lemma,  Propositions
\ref{henselpolyequation} and \ref{henselexpequation} generalize to
their Implicit Function Theorem versions for  {\em systems of
generalized polynomial} and {\em exponential equations} in $\Z_p$.
We set a couple of notations.

\begin{notn}
 We fix  $m, n\in \N,  t_1,\dots,t_m\in \N$,
 and for each $1\le i\le m$ a multivariate generalized polynomial equation
 and a multivariate exponential equation
 $$f_i(x_1,\dots,x_n):=\sum_{j=1}^{t_i}a_j^{(i)}x_1^{\alpha_{j}^{(i1)}}\cdots
 x_n^{\alpha_{j}^{(in)}}\quad \mbox{and} \quad
g_i(\alpha_1,\dots,\alpha_n):=\sum_{j=1}^{t_i}a_j^{(i)}(x_{j}^{(i1)})^{\alpha_{1}}\cdots
 (x_{j}^{(in)})^{\alpha_{n}},$$
 where $a_j^{(i)}\in\Z_p, \alpha_{j}^{(i\ell)}\in\EE_p, x_{j}^{(i\ell)}\in \Z_p^\times,
 $ are given.

\begin{itemize}
\item For a {\em system of $m$ generalized polynomial  equations
in $n$ unknowns $x_1,\dots,x_n\in \Z_p^\times$}:
$$\mathbf{f}(\mathbf{x})=(f_1(x_1,\dots,x_n),\dots,f_m(x_1,\dots,x_n)),$$
 we  denote by $J_p\, \mathbf{f}(\mathbf{x})$ its
{\em pseudo-jacobian} matrix modulo $p$:
$$J_p\, \mathbf{f}(\mathbf{x}):= \left(\begin{array}{ccc}
\Delta_p\,f_1(x_1) &\dots &\Delta_p\, f_1(
x_n)\\
\vdots & & \vdots\\
\Delta_p\, f_m( x_1) &\dots & \Delta_p\,f_m(
x_n)\end{array}\right) \in (\Z/p\Z)^{m\times n}.$$ \item For a
{\em system of $m$ exponential equations in $n$ unknowns
$\alpha_1,\dots,\alpha_n\in \EE_p$}:
$$\mathbf{g}(\mathbf{\alpha})=(g_1(\alpha_1,\dots,\alpha_n),
\dots,g_m(\alpha_1,\dots,\alpha_n)),$$ we denote by
$J_p\,\mathbf{g}(\mathbf{\alpha})$ its {\em pseudo-jacobian}
matrix modulo $p$:
$$J_p\,\mathbf{g}(\mathbf{\alpha}):= \left(\begin{array}{ccc}
\Delta_p\, g_1( \alpha_1 )&\dots &\Delta_p\, g_1
(\alpha_n)\\
\vdots & & \vdots\\
\Delta_p\,g_m( \alpha_1 )&\dots & \Delta_p\, g_m(
\alpha_n)\end{array}\right)\in (\Z/p\Z)^{m\times n},$$
\end{itemize}
where $\Delta_p\, f_i(x_\ell)$ and $\Delta_p\,
 g_i(\alpha_\ell)$ are the corresponding generalizations of  Formula
(\ref{pseudo-derivative}):
$$\begin{array}{lcll}
 \Delta_p\,f_i(x_\ell)& := & \sum_ja_{j}^{(i)}\,\alpha_{j2}^{(i\ell)}
 \,x_\ell^{-1}\,x_1^{\alpha_{j}^{(i1)}}
 \cdots x_n^{\alpha_{j}^{(in)}}& \pmod p \\[3mm]
\Delta_p\, g_i( \alpha_\ell)&:=&
\sum_ja_{j}^{(i)}\,e_2(x_{j2}^{(i\ell)})\,
(x_{j}^{(i1)})^{\alpha_{1}}\cdots (x_{j}^{(in)})^{\alpha_{n}}&
\pmod p.
\end{array}
$$
\end{notn}

\begin{prop} \label{henselsyspolyequation}
 For $m\le n$, let $\mathbf{f}(\mathbf{x})$ denote a system of $m$ generalized polynomial
 equations in $\Z_p$ in $n$ unknowns $x_1,\dots,x_n\in\Z_p^\times$, and let
$\mathbf{y}:=(y_1,\dots,y_m)\in \Z_p^m$. Let
$\mathbf{x}_1=(x_{11},\dots, x_{n1})\in \Z^n$, $p\ndiv x_{\ell
1}$, $1\le \ell\le n$, be such that
\begin{itemize}
\item $\mathbf{f}(\mathbf{x}_1)\equiv \mathbf{y} \pmod p$,
\item $ J_p\, \mathbf{f}(\mathbf{x}_1) \pmod p$ has rank $m$,
\end{itemize}
then there exists  $\mathbf{x}=(x_1,\dots ,x_n )\in
(\Z_p^\times)^n$ such that $x_\ell=\overline{(x_{\ell 1},\dots)}$,
$1\le \ell\le n$, and $\mathbf{f}(\mathbf{x})=\mathbf{y}$ in
$\Z_p$.
\end{prop}

\begin{prop} \label{henselsysexpequation}
For $m\le n$, let $\mathbf{g}(\mathbf{\alpha})$ denote a system of
$m$ exponential
 equations in $\Z_p$ in $n$ unknowns $\alpha_1,\dots,\alpha_n\in\EE_p$,
 and let
$\mathbf{y}:=(y_1,\dots,y_m)\in \Z_p^m$. Let
$\mathbf{\alpha}_1=(\alpha_{11},\dots, \alpha_{n1})\in \Z^n$, be
such that
\begin{itemize}
\item $\mathbf{g}(\mathbf{\alpha}_1)\equiv \mathbf{y} \pmod p$,
\item $J_p\, \mathbf{g}(\mathbf{\alpha}_1) \pmod p$ has rank $m$,
\end{itemize}
then there exists  $\mathbf{\alpha}=(\alpha_1,\dots ,\alpha_n )\in
\EE_p^n$ such that $\alpha_\ell=\overline{(\alpha_{\ell
1},\dots)}$, $1\le \ell \le n$, and
$\mathbf{g}(\mathbf{\bf\alpha})=\mathbf{y}$ in $\Z_p$.
\end{prop}

We prove here Proposition \ref{henselsyspolyequation} since the other
one has exactly the same proof as Proposition \ref{henselexpequation}.
\begin{proof}(Proposition \ref{henselsyspolyequation})
We set
$\mathbf{f}(\mathbf{x})=(f_1(\mathbf{x}),\dots,f_m(\mathbf{x}))$
where
$$f_i(x_1,\dots,x_n):=\sum_{j=1}^{t_i}a_j^{(i)}x_1^{\alpha_{j}^{(i1)}}\cdots
 x_n^{\alpha_{j}^{(in)}}, \ 1\le i\le m.$$ Following the proof of
 Proposition \ref{henselpolyequation}, we assume that
 $\mathbf{x}_k:=(x_{1k},\dots,x_{nk})\in \Z^n$ is constructed such
 that $\mathbf{f}(\mathbf{x}_k)\equiv \mathbf{y}_k \pmod {p^k}$ with
 $\mathbf{x}_k\equiv \mathbf{x}_1 \pmod p$ coordinate-wise. We need
 $x_{\ell(k+1)}= x_{\ell k} + p^k\xi_\ell$ with $\xi_\ell$ to be
 determined for $1\le \ell\le n$ such that
 $\mathbf{f}(\mathbf{x}_{k+1})\equiv \mathbf{y}_{k+1} \pmod
 {p^{k+1}}$. By assumption, for $1\le i\le m$,
$$y_{i(k+1)}-\sum_{j=1}^{t_i}a_{j(k+1)}^{(i)}x_{1k}^{\alpha_{j(k+1)}^{(i1)}}\cdots
 x_{nk}^{\alpha_{j(k+1)}^{(in)}}\equiv p^kz_i \pmod{p^{k+1}}.$$
 Now,
 $$\begin{array}{crclc}
{y}_{i(k+1)}&\equiv & \sum_{j=1}^{t_i} a_{j(k+1)}^{(i)}
 \prod_{\ell=1}^n(x_{\ell k}+p^k
 \xi_\ell)^{\alpha_{j(k+1)}^{(i\ell)}}
& \pmod{p^{k+1}}\\[2mm]
&\equiv&\sum_{j=1}^{t_i}  a_{j(k+1)}^{(i)}
\prod_{\ell=1}^n( x_{\ell
k}^{\alpha_{j(k+1)}^{(i\ell)}}+p^k\,{\alpha_{j(k+1)}^{(i\ell)}}\,
x_{\ell k}^{\alpha_{j(k+1)}^{(i\ell)}-1}\,\xi_\ell)
&\pmod{p^{k+1}}\\[2mm]
&\equiv&\sum_{j=1}^{t_i}  a_{j(k+1)}^{(i)} \big[ \prod_{\ell=1}^n
x_{\ell k}^{\alpha_{j(k+1)}^{(i\ell)}}&
\\
 & & +p^k\,\sum_{\ell=1}^n {\alpha_{j(k+1)}^{(i\ell)}}\, (
\prod_{\ell'\ne \ell} x_{\ell' k}^{\alpha_{j(k+1)}^{(i\ell')}})
x_{\ell k}^{\alpha_{j(k+1)}^{(i\ell)}-1}\,\xi_\ell \ \big]
&\pmod{p^{k+1}}.

\end{array}$$
Hence we need to solve
$$\begin{array}{crclc}
&p^k\,z_i &\equiv & p^k\,\sum_{\ell=1}^n \big[
\sum_{j=1}^{t_i}
 a_{j(k+1)}^{(i)} {\alpha_{j(k+1)}^{(i\ell)}}\, (
\prod_{\ell'\ne \ell} x_{\ell' k}^{\alpha_{j(k+1)}^{(i\ell')}})
x_{\ell k}^{\alpha_{j(k+1)}^{(i\ell)}-1}\big]\xi_\ell
&\pmod{p^{k+1}}\\[2mm]
\iff&z_i &\equiv & \sum_{\ell=1}^n \big[\sum_{j=1}^{t_i}
 a_{j(k+1)}^{(i)} {\alpha_{j(k+1)}^{(i\ell)}}\, (
\prod_{\ell'\ne \ell} x_{\ell' k}^{\alpha_{j(k+1)}^{(i\ell')}})
x_{\ell k}^{\alpha_{j(k+1)}^{(i\ell)}-1}\big]\xi_\ell &\pmod p
\\[2mm]
\iff&z_i &\equiv & \sum_{\ell=1}^n \Delta_p \,f_i({x}_{\ell
1})\,\xi_\ell &\pmod p.
\end{array}
 $$

 Therefore, since by hypothesis, $J_p(\mathbf{f}(\mathbf{x}_1))=
 (\Delta_p \,f_i(x_{\ell 1}))_{i \ell}$ has maximal rank $m$, the
 system has a solution $(\xi_1,\dots,\xi_n)$ modulo $p$.
\end{proof}

Proposition \ref{henselsysexpequation} yields immediately, as a
particular case, another proof of Theorem \ref{zphensel}:

\begin{proof}(Second Proof of Theorem \ref{zphensel})
Let $b\in \Z_p^\times$ be a good basis for taking logarithms,
 $a_j:=b^{\beta_j}$, and consider the exponential system of $2t$
 equations in $2t$ unknowns given by:
 $\mathbf{g}(\beta,\alpha)=(f_1(\beta,\alpha),\dots,
 f_{2t}(\beta,\alpha))$ where $$f_i(\beta,\alpha)=\sum_{j=1}^t
 b^{\beta_j}x_i^{\alpha_j}, 1\le i\le 2t.$$

By Proposition \ref{henselsysexpequation}, this exponential square
system has a solution in $\EE_p$ if $\Rank
J_p\,\mathbf{g}(\beta,\alpha))=2t$, or equivalently that $\det (J_p\,
\mathbf{g}(\beta,\alpha))\not\equiv 0\pmod p$. Furthermore, on this
hypothesis the solution is unique. Now,
$$\Delta_p\,{f_i}( \beta_\ell)= e_2(b_2)b^{\beta_\ell}x_i^{\alpha_\ell}=a_\ell e_2(b_2)
x_i^{\alpha_\ell} \quad \mbox{and} \quad \Delta_p\, f_i(
\alpha_\ell)= e_2(x_{i2})b^{\beta_\ell}x_i^{\alpha_\ell}=
a_\ell\,e_2(x_{i2})x_i^{\alpha_\ell}.$$ Therefore, since $b_2$ is a
primitive root mod$p^2$, $e_2(b_2)\not\equiv 0\pmod{p}$,

$$\det  \left( \begin{array}{cccccc}a_1e_2(b_2)x_1^{\alpha_1} &
\dots & a_te_2(b_2)x_1^{\alpha_t}&a_1e_2(x_{12})x_1^{\alpha_1}&
\dots &
a_t\,e_2(x_{12})x_1^{\alpha_t}\\
 \vdots  & & \vdots  & \vdots & & \vdots \\
a_1e_2(b_2)x_{2t}^{\alpha_1} &
\dots & a_te_2(b_2)x_{2t}^{\alpha_t}&a_1e_2(x_{(2t)2})x_{2t}^{\alpha_1}&
\dots &
a_t\,e_2(x_{(2t)2})x_{2t}^{\alpha_t}
\end{array}\right) \ \equiv$$
$$\equiv a_1\cdots a_t\,e_2(b_2)^t\det
\left( \begin{array}{cccccc}x_1^{\alpha_1} & \dots &
x_1^{\alpha_t}&a_1\,e_2(x_{12})x_1^{\alpha_1}& \dots &
a_t\,e_2(x_{12})x_1^{\alpha_t}\\
 \vdots  & & \vdots  & \vdots & & \vdots \\
x_{2t}^{\alpha_1} & \dots &
x_{2t}^{\alpha_t}&a_1\,e_2(x_{(2t)2})x_{2t}^{\alpha_1}& \dots &
a_t\,e_2(x_{(2t)2})x_{2t}^{\alpha_t}
\end{array}\right)\not\equiv 0\pmod p
$$
by hypothesis.
\end{proof}

\section{Sparse interpolation}
\label{Fewnomial interpolation}

 A polynomial
$f=\sum_{j=0}^da_jx^j\in A[x]$, with $A$ an arbitrary ring, is
usually called a {\em sparse  polynomial} or a {\em fewnomial} if
we focus on its number of non-zero terms, i.e. the number of
$j$'s s.t. $a_j\ne 0$. If it has at most $t$ non-zero terms, i.e.
$f=\sum_{j=1}^t a_jx^{\alpha_j}$, $f$ is called   a {\em
$t$-sparse polynomial} . Here we choose to name such a fewnomial a
{\em $t$-nomial}
 to avoid confusion with the other
usual notions of sparsity. Also we refer to a polynomial with
exactly $t$ non-zero terms as an {\em exact $t$-nomial}.

\smallskip
As mentioned in the introduction, any univariate $t$-nomial in
$\C[x]$ is uniquely determined by its value in $2t$ different
positive values in $\R$, and in \cite{BeTi}, M. Ben-Or and P. Tiwari
produced a beautiful deterministic algorithm that recovers such a
$t$-nomial $f\in \C[x]$ from its value in the $2t$ interpolation
points $x_1:= 1$, $x_2:=a$, $x_3:=a^2$, ..., $x_{2t}:=a^{2t-1},$
where $a$ is not a root of unity of small order.

Furthermore, their algorithm works  for a $n$-multivariate
$t$-nomial $f\in \C[x_1,\dots,x_n]$, using as input interpolation
points
$$\mathbf{x}_1:= (1,\dots,1),\,
\mathbf{x}_2:=(p_1,\dots,p_n),\,\mathbf{x}_3:=(p_1^2,
\dots,p_n^2),\dots,\mathbf{x}_{2t}:=(p_1^{2t-1},\dots,p_n^{2t-1}),$$
where $p_1,\dots,p_n$ are different integer prime numbers.
 The number of arithmetic operations
in $\C$ it performs equals $t^2(\log t+\log(nd))$ where $\log$
denotes the base 2 logarithm. In case the $t$-nomial has integer
coefficients, the bit size of the algorithm is also polynomial in
the maximal bit size $h$ of the coefficients of $f$.

A lot of work has been done in sparse polynomial and rational
interpolation, in different bases of monomials and models. Let us
mention again here some of these, mostly for polynomial
interpolation in the standard monomial basis and the black box
model, meaning you are allowed to choose your interpolation starting
points: the work of M. Ben-Or and P. Tiwari in \cite{BeTi}, of R.E.
Zippel in \cite{Zip90}, of A. Borodin and P. Tiwari in \cite{BoTi},
the series of papers of D. Grigoriev, M. Karpinski and M. Singer in
\cite{GrKaSi90} on finite fields, as well as \cite{ClDrGrKa91}, the
conceptually unifying paper \cite{GKS91} following \cite{DrGr91},
and \cite{GKS94} on rational interpolation. Also
\cite{KaLa88,KLW90,KLL00,Lee01,KaLe03} that improve Ben-Or-Tiwari
algorithms in different ways (producing probabilistic algorithm that
reduce the explosion of size of intermediate integers or  get rid of
the assumption that a bound $t$ for the number of terms is known).

\smallskip
First let us mention  a result that holds when the black-box
representation of the polynomial is in fact a straight-line program,
that is  a program of evaluation that allows constants, and the $+$,
$-$ and $*$ operations. This result does not seem to appear in the
literature although it should be naturally contained in
\cite{GKS91}.  The proof we present here is the fruit of a
discussion with  Michael Singer.  It  is based on the fact  that
 taking derivatives linearizes   exponents,  the same effect than the $e_p$
 ``logarithmic" map we used in the $p$-adic context, and works for
 straight-line programs because of
  Baur and Strassen derivative inequality \cite[(7.7)]{BCS97}.

\begin{teo}\label{slp}
 There is a deterministic algorithm that takes as
input a straight-line program of length $L$ representing a
$t$-nomial $f\in \C[x]$  and returns (the monomial expansion of)
$f$. The complexity of the algorithm is of order $O(t^4L)$.
\end{teo}

\begin{proof} For $f=\sum_{j=1}^ta_jx^{\alpha_j}\in\C[x]$ we set
$Df(x):=xf'(x)$ so that $D(x^\alpha)=\alpha x^\alpha$, the
$k$-iteration $D^{(k)}(x^\alpha)=\alpha_j^kx^\alpha$, and finally
for $k\in\N$, $D^{(k)}f=\sum_j a_j\alpha^kx^{\alpha_j}$.

Since $f$ is represented by a straight-line program of length $L$,
$\{D^{(k)}f,0\le k\le t\}$ are given by $t+1$ straight-line programs
of length $O(L)$, that can all be constructed from the straight-line
program for $f$ in time $O(tL)$ \cite[(7.7)]{BCS97}.

We fix different positive $x_1,\dots,x_t\in \R$ and we construct the
matrix

$$C(f):=\left(\begin{array}{ccc}
f(x_1)&\dots &f(x_t)\\
Df(x_1)&\dots &Df(x_t)\\
\vdots & & \vdots\\
D^{(t-1)}f(x_1)&\dots &D^{(t-1)}f(x_t)
\end{array}\right)\in \C^{t\times t}$$
and observe that $$C(f)=V(f)\,A(f)\,W(f),$$
where
$$V(f):=\left(\begin{array}{ccc}\alpha_1^0&\dots & \alpha_t^0\\
\vdots & & \vdots\\
\alpha_1^{t-1}& \dots & \alpha_t^{t-1}\end{array}\right),\
A(f):=\left(\begin{array}{ccc}a_1& &\\
& \ddots& \\ & & a_t\end{array}\right),
\ W(f):=\left(\begin{array}{ccc}x_1^{\alpha_1}&\dots & x_t^{\alpha_1}\\
\vdots & & \vdots\\
x_1^{\alpha_t}& \dots & x_t^{\alpha_t}\end{array}\right).
$$

Therefore the rank of $C(f)$ gives the exact number of non-zero coefficients
$a_j$ and we can assume $f$ is an exact $t$-nomial, so that $C(f)$ is invertible.

Now we observe that if we set
$g:=Df=\sum_{j=1}^ta_j\,\alpha_j\,x^{\alpha_j}$, then
$C(g)=V(g)A(g)W(g)$ where $V(g)=V(f),\ W(g)=W(f)$ and $A(g)$ is the
diagonal matrix with diagonal terms $a_j\alpha_j$. Therefore
$$C(g)C(f)^{-1}=V(f)\left(\begin{array}{ccc}
\alpha_1&  & \\
& \ddots &\\
& & \alpha_t\end{array}\right)V(f)^{-1},$$ and to compute the
exponents $\alpha_j,\ 1\le j\le t$, is is enough to compute the
characteristic polynomial of $C(g)C(f)^{-1}$ and its (integer)
roots. The coefficients are then recovered by solving a linear
system.
\end{proof}

\bigskip
In next section  we  investigate  univariate fewnomials over
finite fields, in particular in $(\Z/p\Z)[x]$, and over the finite
rings $\Z/p^k\Z$, where $p$ is an odd prime number.

\smallskip
We then switch to univariate integer polynomials, trying to produce
an answer to the question raised by Ben-Or and  Tiwari on how to
interpolate a $t$-nomial in $\Z[x]$ from $2t$ arbitrary different
real positive values.

\smallskip
In the sequel we denote by
 $\deg f$
  the degree of a non-zero polynomial $f\in A[x]$, $A$ a ring,  and by $h(f)$ its binary
length when $A=\Z$, i.e. the maximum (base $2$) logarithm $\log$
of the absolute value of its coefficients.

\subsection{Fewnomial interpolation in $(\Z/p\Z)[x]$ and $(\Z/p^k\Z)[x]$}

{\ }\\
Here we  impose conditions on the exponents of the
 polynomial, since for example, every $x\in (\Z/p\Z)[x]$ is a root of the binomial
 $f=x^p-x$.
We first recall the proof of a well-known   result (see for instance
\cite[Th 4.2]{ClDrGrKa91} for a more general context): any
$t$-nomial of degree strictly bounded by $p-1$ in $(\Z/p\Z)[x]$ is
uniquely determined by its values in $2t$ points of the form
$x_i=\rho^{i-1}$ for $1\le i\le 2t$, where $\rho$ is a primitive
root modulo $p$. (Of course this statement makes sense only if $t$
is small with respect to $p$, i.e. if $2t\le p-2$, since if not
standard interpolation suffices.) In some sense this is an analogous
statement than that of a $t$-nomial in $\C[x]$ being uniquely
determined by its value in $2t$ positive points. This condition can
not be arbitrarily relaxed since for example  a binomial modulo $7$
of degree $< 6$ is not uniquely determined by its values in $1, 2,
3$ and $4$: $f=x^4+3$ and $g=3x^3+x$ coincide modulo $7$ but are
different (observe that in this example we may take $\rho=3$, and
therefore $1$, $3$, $2$ and $6$ are good interpolation points:
$f(6)\ne g(6)$).

\smallskip
\begin{obs} Let $p$ be a prime number and $\rho\in \Z$ be a primitive root modulo $p$.
Let $f=\sum_{j=1}^t a_jx^{\alpha_j} \in (\Z/p\Z)[x]$  be a
$t$-nomial satisfying that for $j\ne \ell$, $\alpha_j\not \equiv
\alpha_\ell\pmod{(p-1)}$. Then $f(\rho^i)= 0$ for $0\le i\le t-1$
implies $f= 0$.
\end{obs}

\begin{proof}
We have
$$\left(\begin{array}{ccc}1^{\alpha_1}& \dots & 1^{\alpha_t}\\
\rho^{\alpha_1}& \dots & \rho^{\alpha_t}\\
\vdots &  & \vdots\\
\rho^{(t-1)\alpha_1}& \dots & \rho^{(t-1)\alpha_t}
\end{array}\right) \,\left(\begin{array}{c}a_1\\ \vdots \\ a_t
\end{array}\right)=\left(\begin{array}{c}0\\\vdots \\ 0
\end{array}\right).$$
But the Vandermonde determinant of the left-hand side matrix
$$\prod_{1\le j<\ell\le t} (\rho^{\alpha_\ell}-\rho^{\alpha_j})$$

does not vanish modulo $p$ since $\rho$ is a primitive root modulo
$p$ and $p-1\ndiv \alpha_\ell-\alpha_j$. Thus the unique solution
of the system is $a_j= 0$ for $1\le j\le t$.
\end{proof}

\begin{cor} \label{p-identidad}
 Set $p$ a prime number. A $t$-nomial in $(\Z/p\Z)[x]$
 of degree strictly bounded by $p-1$ is
uniquely determined  by its values at $\rho^i$, $0\le i\le 2t-1$,
where $\rho$ is a primitive root modulo $p$.
\end{cor}

\begin{proof} Let $f:=\sum_{j=0}^{p-2}a_jx^{j}, g:=
\sum_{j=0}^{p-2}b_jx^{j}$ be  two $t$-nomials in $(\Z/p\Z)[x]$
such that $f(\rho^i)=g(\rho^i)$, $0\le i\le 2t-1$. Then $h:=
\sum_{j=0}^{p-2}(a_j-b_j)x^j$  is a $2t$-nomial modulo $p$ that
satisfies $h(\rho^i)\equiv 0\pmod p$ for $0\le i\le 2t-1$, and we
apply the previous observation.
\end{proof}

In the  paper mentioned above, M. Ben-Or and P. Tiwari raised the
problem of generalizing their procedure for finite fields. As a
partial answer, it is straight-forward that the determination of the
unique  $t$-nomial  in $(\Z/p\Z)[x]$ of degree strictly bounded by
$p-1$ and  with prescribed values in $1, \rho, \dots, \rho^{2t-1}$
---when it exists--- can be easily done copying their algorithm for
this case. This gives an alternative simple proof in this
particular case of the multivariate result of \cite[Th.
4.2]{ClDrGrKa91}.

\begin{alg}  \label{BenOr} {\em (Ben-Or/Tiwari)}

Set $p$ a prime number and $\rho \in \Z$ a primitive root modulo
$p$. Let $f\in (\Z/p\Z)[x]$ be a $t$-nomial of degree strictly
bounded by $p-1$.
 Then there is a deterministic
algorithm that takes as inputs the values $y_1:=f(1),y_2:=f(\rho),
\dots,y_{2t}=f(\rho^{2t-1})$  and returns  $f$.  The binary
running time of the algorithm equals $O((t^2+p)\log p)$
\end{alg}

\begin{proof} The proof copies \cite{BeTi}.
Let $f=\sum_{j=1}^t a_jx^{\alpha_j}\in (\Z/p\Z)[x]$, $0\le
\alpha_j\le p-2$. The algorithm first computes the exact number
$\tilde t$ of terms of $f$, then it determines the exponents
$\alpha_1,\dots,\alpha_{\tilde t}$ associated to non-zero
coefficients and finally it recovers the coefficients
$a_1,\dots,a_{\tilde t}$.

\smallskip
Let us first assume that $\tilde t=t$. The core of the procedure
is  the same previous fact that  since $\rho $ is a primitive root
modulo $p$,  $\rho^{\alpha}\ne \rho^\beta $ in $\Z/p\Z$ for $0\le
\alpha\ne \beta \le p-2$.

As in \cite{BeTi}, we construct a polynomial $F\in
(\Z/p\Z)[\lambda]$ of degree $t$ whose roots are exactly
$\rho^{\alpha_j}$, and then we recover $\alpha_j$, $1\le j\le t$,
by simple inspection.

The polynomial
$F=\prod_{j=1}^t(\lambda-\rho^{\alpha_j})=\sum_{k=0}^t
b_k\lambda^k$, $b_t=1$, is constructed in the following way: For
$0\le \ell\le t-1$, $1\le j\le t$: \begin{eqnarray*} 0&=&
a_j\rho^{\alpha_j\ell}F(\rho^{\alpha_j})=a_j(b_0\rho^{\alpha_j\ell}+
b_1\rho^{\alpha_j(\ell+1)} + \cdots +
b_t\rho^{\alpha_j(\ell+t)})\quad \Longrightarrow \\
0&=& \sum_{j=1}^ta_j\rho^{\alpha_j\ell}F(\rho^{\alpha_j})\\
&=& b_0\sum_{j=1}^t a_j\rho^{\alpha_j\ell} + b_1\sum_{j=1}^t
a_j\rho^{\alpha_j(\ell+1)}+\cdots + b_t\sum_{j=1}^t
a_j\rho^{\alpha_j(\ell+t)}\\
&=& b_0 f(\rho^\ell) + b_1 f(\rho^{\ell+1}) + \cdots +
b_{t-1}f(\rho^{\ell+t-1}) + f(\rho^{\ell+t})\\[3mm]
&=& b_0 y_{\ell+1} +\cdots + b_{t-1}y_{\ell+t} + y_{\ell+t+1}.
\end{eqnarray*}
This yields the following system
$$\left(\begin{array}{ccc}y_1&\dots&y_t\\
\vdots& & \vdots\\
y_t&\dots & y_{2t-1}
\end{array}\right)\,\left(\begin{array}{c}b_0\\\vdots \\
b_{t-1}
\end{array}\right)=- \left(\begin{array}{c}y_{t+1}\\
\vdots\\ y_{2t}
\end{array}\right)$$
which is clearly solvable since
$$\left(\begin{array}{ccc}y_1&\dots&y_t\\
\vdots& & \vdots\\
y_t&\dots & y_{2t-1}
\end{array}\right)=\left(\begin{array}{ccc}1&\dots&1\\
\rho^{\alpha_1}&\dots&\rho^{\alpha_t}\\
\vdots& & \vdots \\
\rho^{\alpha_1(t-1)}&\dots &\rho^{\alpha_t(t-1)}
\end{array}\right) \left(\begin{array}{ccc}a_1& & \\
 &\ddots &  \\
 &  & a_t
\end{array}\right) \left(\begin{array}{ccc}
1&
\dots&\rho^{\alpha_1(t-1)}\\
1&
\dots&\rho^{\alpha_2(t-1)}\\
\vdots& & \vdots\\
1&
\dots&\rho^{\alpha_t(t-1)}
\end{array}\right)$$
whose determinant $a_1\dots a_t \prod_{1\le i<j\le
t}(\rho^{\alpha_j}-\rho^{\alpha_i})^2\not\equiv 0\pmod p$.

Now it is easy to recover $a_1,\dots,a_t$ by solving the
 Vandermonde system:
$$\left(\begin{array}{ccc}1&\dots&1\\
\rho^{\alpha_1}&\dots &\rho^{\alpha_t}\\
\vdots & & \vdots\\
\rho^{\alpha_1(t-1)}&\dots &
\rho^{\alpha_t(t-1)}\end{array}\right) \left(\begin{array}{c}
a_1\\\vdots \\a_t\end{array}\right)=\left(\begin{array}{c}
y_1\\\vdots \\y_{t}\end{array}\right)$$

\smallskip
When the exact number of terms $\tilde t$ is not known, but
strictly bounded by $t$, we compute it by considering for $1\le
\ell\le t$ the matrices $$V_\ell:=\left(\begin{array}{ccc}
y_1&\dots & y_\ell\\\vdots & & \vdots \\y_{\ell}&\dots &
y_{2\ell-1}\end{array}\right),$$ that satisfy that
$\det(V_\ell)=0$ for all $\tilde t<\ell\le t$ while
$\det(V_{\tilde t})\ne 0$ in $\Z/p\Z$. Thus, computing $\tilde t$
is equivalent here to compute the rank of the matrix $V_t$.

\medskip
Let us check the complexity of this
 algorithm: all integers are bounded by $p$, and $t<p$.
 As in \cite{BeTi}, computing the exact number $\tilde t$ of terms of the $t$-nomial
 requires
 $O(t^2\log p)$ binary operations, and bounds the complexity of
 solving
the linear system to determine the polynomial $F$. The simplest
way of  computing  the  exponents $\alpha_j$,
 $1\le j\le t$, of $f$ seems to be by simple inspection: computing
 $\rho^{i}$, $1\le i\le p-2$ and checking which of those are roots
 of $F$.
This takes $O(p\log p)$ steps. Finally, recovering the
coefficients does not modify the overall complexity.
\end{proof}

\begin{cor} An analogous algorithm  holds in a finite field $\F_q$ for $q=p^n$, since the
multiplicative group of a finite field is cyclic of order $q-1$.  Any
$t$-nomial in $\F_q[x]$ of degree strictly bounded by $q-1$ can be
recovered from its values in the interpolation points $\rho^{i}$,
$0\le i\le 2t-1$, where $\rho$ is a generator of the field $\F_q$ over
$\Z/p\Z$.
\end{cor}

Now we turn to polynomials with coefficients in the ring $\Z/p^k\Z$,
where all usual arguments fail for it is not even a domain. However
combining Theorem  \ref{hensel} and Algorithm \ref{BenOr}, we are
able to obtain some results for $p$ an odd prime number and $k\in
\N$, $k\ge 2$. We introduce for polynomials in $(\Z/p^k\Z)[x]$ the
analogue of Definition \ref{reduceswell}:

\begin{defn}\label{reduceswellpk} Let $p$ be an odd prime number and $k\in
\N$. We say that a polynomial $f=\sum_j a_jx^{\alpha_j}\in
(\Z/p^k\Z)[x]$ with $a_j\ne 0$, $\forall \, j$, {\em reduces well
modulo} $p$ if $p\ndiv a_j$ for any $j$,  and $p-1\ndiv
\alpha_j-\alpha_\ell$ for any $j\ne \ell$.
\end{defn}

\begin{cor} \label{pkidentidad} Set $p$ an odd prime number and $k\in \N$ with $k\ge 2$.
A $t$-nomial in $(\Z/p^k\Z)[x]$ of degree strictly bounded by
$\varphi(p^k)$ that reduces well modulo $p$ is uniquely determined
by its values at $\rho^i$, $0\le i\le 2t-1$, where $\rho$ is a
primitive root mod$p^2$.
\end{cor}

\begin{proof}
Let $f=\sum_{j=1}^t a_jx^{\alpha_j+(p-1)k_j} \in (\Z/p^k\Z)[x]$ be
such a $t$-nomial, where $a_j\not\equiv 0\pmod p$ and $0\le
\alpha_j< p-1$ are all distinct since  $f$ reduces well modulo $p$.
Since $\rho$ is also a primitive root modulo $p$, by Corollary
\ref{p-identidad}, $\sum_j a_jx^{\alpha_j}\in (\Z/p\Z)[x]$ is the
unique (exact) $t$-nomial  of degree bounded by $p-1$ with the
prescribed values in $\rho^i$, $0\le i\le 2t-1$. Applying Theorem
\ref{hensel}, since $\{\rho^{i}, 0\le i\le 2t-1\}$ is a good
starting set (Proposition \ref{startingset}), there exists a unique
$g$ such that $g(\rho^{i})\equiv f(\rho^i)\pmod {p^k}$, $0\le i\le
2t-1$, under the condition that the coefficients coincide modulo
$p^k$ and the exponents mod$\varphi(p^k)$. Therefore, since $f$ is
such a polynomial,  $f\in (\Z/p^k\Z)[x]$ is the unique $t$-nomial of
degree bounded by $\varphi(p^k)$ that reduces well modulo $p$.
\end{proof}

\begin{alg}  \label{pkBenOr}

Set $p$  an odd prime number, $\rho \in \Z$  a primitive root
mod$p^2$ and $k\in \N$ with $k\ge 2$. Let $f\in (\Z/p^k\Z)[x]$ be a
$t$-nomial of degree strictly bounded by $\varphi(p^k)$ that reduces
well modulo $p$.
 Then there is a deterministic
algorithm that takes as inputs the values $y_1:=f(1),y_2:=f(\rho),
\dots,y_{2t}=f(\rho^{2t-1})$ in $\Z/p^k\Z$  and returns  $f$.  The
binary running time of the algorithm is of order  $O(t^3k^2 \log^2
p+p\log p)$.
\end{alg}

\begin{proof}
We first apply Algorithm \ref{BenOr} to compute the unique exact
$t$-nomial $f_0\in(\Z/p\Z)[x]$ of degree strictly bounded by $p-1$
such that $f_0(\rho^{i-1})\equiv y_i\pmod p$ for $0\le i\le 2t-1$.
Then we apply Theorem \ref{hensel} to lift $f_0$ to $f$.

Now we check the complexity of the algorithm. First we reduce $y_i$
and $\rho^{i-1}\pmod p$, $1\le i\le 2t$ to construct $f_0$. This
takes $O(t\,k\log p + (t^2+p)\log p)$ binary operations. Next, there
are at most $\log\log \varphi(p^k)=O(\log k + \log \log p)$ lifting
steps, each of them solves a system of size $2t$ with entries
bounded by $p^k$, that takes $O((\log k+\log \log p)t^3 k\log p)$
binary operations. The overall complexity is then of order
$O(t^3k\log k \log p+ t^3k\log p\log \log p+p\log p)$.
\end{proof}

 A final observation for this section is that, with the same
proof than  that of Corollary \ref{pkidentidad}, Proposition
\ref{pkequivalencia} can be reformulated as follows:

\begin{cor} Set $p$  an odd prime number.
Let $f=\sum_j a_jx^{\alpha_j},g=\sum_\ell b_\ell
x^{\beta_\ell}\in\Z[x]$ be two polynomials that reduce well
modulo $p$.
 Then, for any $k\in \N$,
 the three following conditions are equivalent:
\begin{itemize}
\item $f$ and $g$ have the same number $t$ of non-zero terms, and
up to an index permutation, $a_j\equiv b_j\pmod{p^k}$ and
$\alpha_j\equiv \beta_j \pmod{\varphi(p^k)}$.
 \item $f(x)\equiv g(x) \pmod{p^{k}}$
for all
 $x\in \Z$ prime to $p$.
 \item $f(\rho^{i-1})\equiv g(\rho^{i-1})\pmod {p^k}$ for   $1\le i\le
 2t$ and
 $\rho\in \Z$ a primitive root modulo $p^2$.
\end{itemize}
\end{cor}

For the previous results we used the fact that for $\rho\in \Z$ a
primitive root modulo $p^2$,  $\{1,\rho,\dots,\rho^{2t-1}\}$ is a
good starting set (of type (1) of Proposition \ref{startingset}).
For good starting sets of type (2), for instance
$\{1,\rho,\dots,\rho^{t-1}, p+1,p+\rho, \dots,p+\rho^{t-1}\}$ for
$\rho\in \Z$ a primitive root modulo $p$, we do not have  analogous
of Corollary \ref{p-identidad} and Algorithm \ref{BenOr}, and the
best we can obtain are the following statements:

\begin{cor}  Set $p$ an odd prime number and $k\in \N$, $k\ge 2$.
A $t$-nomial in $(\Z/p^k\Z)[x]$ of degree strictly bounded by
$\varphi(p^k)$ that reduces well modulo $p$ is uniquely determined by
its values in $\{1,\dots,p-1,p+1,\dots,2p-1\}$.
\end{cor}

\begin{proof}
This is simply due to the fact that the first $p-1$ points decide
which is the $t$-nomial modulo $p$ and then we apply Theorem
\ref{hensel} for the  over-constrained compatible system we have,
using that the given set contains a good starting set of type
 (2).
\end{proof}

\begin{alg}
Set $p$  an odd prime number and  $k\in \N$, $k\ge 2$. Let $f\in
(\Z/p^k\Z)[x]$ be a $t$-nomial of degree strictly bounded by
$\varphi(p^k)$ that reduces well modulo $p$.
 Then there is a deterministic
algorithm that takes as inputs the values $f(1),\dots, f(p-1),
f(p+1), \dots, f(2p-1)$ in $\Z/p^k\Z$ and returns $f$. The binary
running time of the algorithm is of order  $O(p^4k\log k )$.
\end{alg}

\begin{cor} Set $p$  an odd prime number.
Let $f=\sum_j a_jx^{\alpha_j},g=\sum_\ell b_\ell
x^{\beta_\ell}\in\Z[x]$ be two polynomials that reduce well
modulo $p$.
 Then, for any $k\in \N$,
 the three following conditions are equivalent:
\begin{itemize}
\item $f$ and $g$ have the same number $t$ of non-zero terms, and
up to an index permutation, $a_j\equiv b_j\pmod{p^k}$ and
$\alpha_j\equiv \beta_j \pmod{\varphi(p^k)}$.
 \item $f(x)\equiv g(x) \pmod{p^{k}}$
for all
 $x\in \Z$ prime to $p$.
 \item $f(i)\equiv g(i)\pmod {p^k}$ for   $1\le i\le
 p-1$ and $p+1\le i\le 2p-1$.
\end{itemize}
\end{cor}

\smallskip
\subsection{Fewnomial interpolation in $\Z[x]$ }

{\ } \\
In their paper, M. Ben-Or and P. Tiwari also  raised the problem of
producing an algorithm that interpolates a $t$-nomial in $\C[x]$
from $2t$ arbitrary different real positive values. Here we restrict
to polynomials in $\Z[x]$. On one hand we observe that applying a
bound by A. Borodin and P. Tiwari \cite[Thm.4.3]{BoTi} we can
restrict ourselves to $t+1$ interpolation  points, $t$ of them being
almost arbitrary, but the last one imposed and huge. On the other
hand, Theorem \ref{hensel} enables us to  reduce the size of the
starting interpolation points for  $t$-nomials in $\Z[x]$
  that reduce well mod $p$ for some small enough prime number $p$.

\begin{teo}\label{Borodin} \   {\em (\cite[Thm.4.3]{BoTi})}

Let  $f\in \Z[x]$ be a $t$-nomial, and $(x_i,y_i)\in \Z^2, \ 1\le
i\le t$,   that satisfy that $x_i\ge 2$ and $f(x_i)=y_i$ for $1\le
i\le t$. Then
$$\deg f  \le \max_i\log x_i + t^2\max_i\log y_i +2.$$
\end{teo}

This bound for the degree of such a $t$-nomial $f$ in term of the
height of its evaluation points immediately yields a bound for the
height $h(f)$ of $f$:

\begin{cor} Let  $f\in \Z[x]$ be a $t$-nomial, and $(x_i,y_i)\in \Z^2, \ 1\le
i\le t$,   that satisfy that $x_i\ge 2$ and $f(x_i)=y_i$ for $1\le
i\le t$. Then
 $$h(f)\le \max_i\log|y_i| + 2t(\log t +\max_i\log x_i)
+ t\,{\max_i}^2\log x_i + t^3\max_i\log x_i\max_i\log|y_i|.$$
\end{cor}

\begin{proof} Set $D:=\max_i\log x_i + t^2\max_i\log y_i +2$.
 As $f=\sum_{j=1}^ta_jx^{\alpha_j}$ where $\alpha_j\le D$, when solving the
  linear system  induced by $f(x_i)=y_i$ for $1\le i\le t$, we deal with  an
  integer matrix of size $t$ and entries of absolute height bounded by
  $H:=t\,\max_i{x_i}^D$ and a vector with entries of absolute height bounded by $\max_i |y_i|$.
  Applying Cramer's rule and the fact that $a_i\in \Z$,
  we obtain $|a_i|\le t! \,H^{t-1}\,\max_i|y_i|$. This gives the announced bound.
\end{proof}

Now we remark that an a priori bound  for the height of a $t$-nomial
$f\in \Z[x]$
  immediately  yields the polynomial by
 interpolation in
one single huge value:

\begin{obs} Let $f\in \Z[x]$ be a  $t$-nomial and let $H$ be a
 bound for the maximum absolute value of the coefficients of $f$.

Then, for any odd number $\widetilde H\ge  2H+1$, there is a
deterministic algorithm that takes as input the value
$f(\widetilde H)$ and returns the $t$-nomial $f$ in $\Z[x]$.

The binary running time of the algorithm is polynomial in
   $t$, $\log(\widetilde H)$ and $\log(f(\widetilde H))$.
\end{obs}

\begin{proof} It is enough to write
$f(\widetilde H)=\sum_{i\ge 1}   \tilde a_i\widetilde H^{i}$ in
base $\widetilde
 H$ with coefficients in $[-\frac{\widetilde H-1}{2},\frac{\widetilde H-1}{2}]$
 to
recover the exponents and coefficients of $f$ by simple
inspection.
\end{proof}

Combining these two  facts,  we conclude that to obtain $f$ we can
restrict ourselves to $t+1$ interpolation  points
$x_1,\dots,x_{t+1}$, where $x_i\ge 2$ are arbitrary for $1\le i\le
t$ but $x_{t+1}$ satisfies the  condition of the previous
observation.

\bigskip
Now we come back to  Ben-Or/Tiwari type algorithms to recover
$t$-nomials. We remind that in the sequel a bound $t$ for the number
of non-zero terms is always given as an input.  An inconvenient in
the original algorithm by \cite{BeTi} is the explosion of
intermediate integers: one has to deal with a polynomial $F$
---see proof of Algorithm \ref{BenOr} above--- where a coefficient
equals
 at least $2^{\alpha_1+\dots + \alpha_t}$. This problem has been  solved
 in \cite{KaLa88, KLW90}, where the authors
propose  a probabilistic algorithm that  keeps the intermediate
integers small by employing the  traditional Hensel lifting of
roots. They choose  a ``lucky" prime $p$ and $k\in \N$ big enough
(essentially s.t. $p^k\ge \deg f $), compute  the crucial polynomial
$F$ modulo $p^k$ and lift  its roots modulo $p$ to roots modulo
$p^k$. These algorithms require a degree bound as input,   mostly to
control the probability of unlucky reduction modulo $p$.

Here we present an alternative algorithm
 to recover a $t$-nomial $f\in \Z[x]$
from its interpolation in $2t$ points of size bounded by $p^2$,
provided we know in advance that  it reduces well modulo the odd
prime number $p$. Since we are still not able to produce a
probability analysis for the choice of a good prime $p$ such that
$f$ reduces well modulo $p$, our algorithm only yields an heuristic
for arbitrary $t$-nomials.

We de not intend here to compare the speed of our method with that
of  \cite{KLW90}: no serious implementation has been done yet.
However, since both methods are different in nature, we think that
it  can be useful to have them both in mind.

\smallskip

\begin{alg}\label{interpolation}

Let $f\in \Z[x]$ be a  $t$-nomial and let $p>t$ be an odd prime
number {such that $f$ reduces well  modulo $p$}. Set $\rho \in \Z$ a
primitive root modulo $p^2$   and let $x_1,\dots,x_{2t}\in \N$ be
such that $x_i\equiv \rho^{i-1}\pmod {p^2}$.

Then there is a deterministic algorithm that takes as input the
values $y_1:=f(x_1),\dots, y_{2t}:=f(x_{2t})$ and returns the
 $t$-nomial $f$ in $\Z[x]$.

The binary running time of the algorithm is polynomial in $p$,
$\log d$, $h$ and $\tilde h$, where $d:=\deg f$, $h:=h(f)$ and
$\tilde h:=\max \{h(y_i), 1\le i\le 2t\}$.

The algorithm  computes  $m\le \max\{\, \lceil \log \log d\rceil
,\lceil \log h\rceil \,\}$   $t$-nomials $f_0,\dots,f_m\in \Z[x]$
until matching $f$. The termination of the procedure is given by the
condition $f_m(x_i)=y_i$ for $1\le i\le 2t$.
\end{alg}

\begin{proof}
We first compute by Algorithm \ref{BenOr} the  unique exact $\tilde
t$-nomial $f_0\in \Z[x]$, where $\tilde t\le t$, of degree $\le p-2$
and integer coefficients in $[-\frac{p-1}{2},
 \frac{p-1}{2}]$,
 determined by the conditions
 $$f_0(x_i)\equiv y_i \pmod p \quad \mbox{for} \ 1\le i \le 2t.$$

 This $\tilde t$-nomial must exist since it coincides with
 $\sum_{j=1}^{t}a_jx^{\alpha_j}$ if
 $f:=\sum_{j=1}^{t}(a_j+p\,d_j)x^{(\alpha_j+(p-1)\delta_j)}$,
 with $a_j\in [-\frac{p-1}{2},
 \frac{p-1}{2}], \beta_j\in[0, p-2]$.

We observe that if $f_0(x_i)=y_i$ in $\Z$, the procedure stops and
$f=f_0$, since $f$ is uniquely determined by its value in the $2t$
positive values $x_1,\dots,x_{2t}$.

\medskip
W.l.o.g. we can assume now that $\tilde t=t$. To continue the
procedure we apply Theorem \ref{hensel} to compute recursively the
unique exact $t$-nomial $f_k=\sum_{j=1}^tb_jx^{\beta_j}$ of degree
strictly bounded by $\varphi(p^{2^k})$ and with integer coefficients
in $[ -\frac{p^{2^k}-1}{2}, \frac{p^{2^k}-1}{2}]$ that satisfies
$f_k(x_i)\equiv y_i\pmod {p^{2^k}}$ for $1\le i\le 2t$.

\smallskip
The termination of the procedure occurs at most for $f_{m}=f$, i.e
$m$ such that
   $\frac{p^{2^{m}}-1}{2}\ge 2^h$
and $\varphi(p^{2^{m}})>d$, that is $$m= \max\{\, \lceil \,\log h
\, \rceil, \lceil \,\log \log d \,\rceil\,\}.$$

\smallskip
Now let us compute the binary complexity of the algorithm.  The
running time needed to compute $ f_0$ is of order $O((t^2+p)\log
p)$.

  To compute $ f_{k+1}$ from $f_k$:

First we compute mod$p^{2^k}$ the entries of matrix $M_k$: we need
to compute $e_{2^{k+1}}(x_i)$ and $\ell_{2^{k+1}}(x_i)$ defined in
Identity (\ref{eielei}), that require
$O(t\,(\log(p^{2^{k+1}})+\tilde h))=O(t(2^{k+1}\log p +\tilde h))$
bit operations. To compute $x_i^{\beta_j}$ modulo $p^{2^k}$ requires
 $O(\log(p^{2^k}))= O(2^k\log p)$ operations.
The computation of the  determinant of $M_k$ and of its inverse
modulo $p^{2^k}$ requires $O(t^32^k\log p)$ more operations. Thus,
computing $f_{k+1}$ from $f_k$ requires $O(t^3(2^{k+1}\log p+\tilde
h))$ bit operations.

Thus, the total number of bit operations of the algorithm is
bounded by
$$O(t^3(\max\{ h, \log d\,\}\log p +\tilde h) + p\log p).$$
(We kept the complexity  in terms of $p$ and $\log p$ since the
only place where it seems to depend on $p$ is in the computation
of the starting polynomial $f_0$).
\end{proof}

 Since under our conditions, the bound of
\cite[Thm.4.3]{BoTi} gives $\deg f\le 2\log p +t^2\tilde h + 2$, we
obtain the following heuristic for the case   we do not know in
advance that $f$ reduces well modulo $p$.

\begin{heur}

{\ }

\begin{itemize}
\item
{\bf Input:}  $f\in \Z[x]$ given by a black-box,
 $t\in \N $ a bound for the number of terms of $f$.
\item{\bf
Output:} Luckily, the monomial basis representation of $f$.
\item{\bf Heuristic:}
\begin{itemize}
\item Pick $p>t$ an odd prime
number.
\item Pick $\rho\in \Z$  a primitive root modulo $p^2$.
\item
Pick $x_1,\dots,x_{2t}\in \N$  such that $x_i\equiv \rho^{i-1}\pmod
{p^2}$.

\item  Compute  $y_i:=f(x_i)$ for $1\le i\le 2t$ from the black box.
\item  Compute $f_0$,
the unique exact $\tilde t$-nomial  modulo $p$ such
that $f_0(x_i)\equiv y_i\pmod p$ for $1\le i\le 2t$ (we observe
that $\tilde t\le t$ must occur).
\item Lift,
Applying Theorem \ref{hensel}, $f_0$ to
$f_m=\sum_{i=1}^{\tilde t} a_ix^{\alpha_i}$ such that
$$\varphi(p^{2^{m-1}})\le 2\log p +t^2\tilde h + 2<\varphi(p^{2^{m}}). $$

(This yields the possible exponents $\alpha_1,\dots, \alpha_{\tilde
t}$ of $f$.)
\item  Set $ \tilde f=\sum_{i=1}^{\tilde t} z_ix^{\alpha_i}$
and try to interpolate $\tilde  f(x_i)=y_i$ for $1\le i\le 2t$ in
$\Z[x]$ solving a simple  Vandermonde system.

\item If the interpolation  problem has
a solution, then $\tilde f=f$ and output $\tilde f$.

\item  If there is no solution, it was because $f$ was
not an exact $\tilde t$-nomial (and in fact the exact number of
terms of $f$ is strictly greater than $\tilde t$). In that case
pick another prime $q>t$ and start
 the procedure again.
(If the exact number of terms of the new starting polynomial $f_0$ is not greater than
 $\tilde t$, pick another prime.)\hfill\mbox{$\Box$}
 \end{itemize}
 \end{itemize}
\end{heur}

\bigskip

\noindent{\bf Final comment:} The problem of finding an algorithm
that, given an (unknown) $t$-nomial $f$ in $\Z[x]$ and $2t$ starting
evaluation points, finds the monomial structure of $f$ has proven
very hard. If one can find an algorithm that solves the (easier)
problem over the finite field $\Z/p\Z$ (where there may be no
solution or more than one), our method for lifting the coefficients
and the exponents can be used under the assumption of the existence
of a ``good'' (relatively small) prime $p$, i.e. a prime $p$ such
that $f$ reduces well modulo $p$, and the pseudo-jacobian of $f$ is
invertible. A probability analysis for good reduction of $t$-nomials
modulo $p$ is still lacking. We are trying to give an answer to
these problems.

\end{document}